\documentclass[authoryear,review]{elsarticle}
\usepackage[margin=2.5cm]{geometry}
\usepackage{subcaption}
\usepackage{lineno}
\usepackage{hyperref}
\usepackage{xurl}
\usepackage{multicol}

\usepackage[noend]{algpseudocode}
\usepackage{amssymb,amstext,amsmath,amsthm,mathrsfs}
\usepackage{graphicx,caption,varwidth,rotating}
\usepackage[linesnumbered,ruled,vlined]{algorithm2e}
\usepackage[table]{xcolor}

\SetKwInput{KwInput}{Input}                
\SetKwInput{KwOutput}{Output}              

\biboptions{authoryear}
\begin{document}

\begin{frontmatter}

\title{A matheuristic for tri-objective binary integer programming}

\author[1]{Duleabom An\corref{cor1}}
\ead{duleabom.an@jku.at}\cortext[cor1]{Corresponding author}

\author[1]{Sophie N. Parragh} \author[1]{Markus Sinnl}%
\author[2]{Fabien Tricoire}
\address[1]{Institute of Production and Logistics Management/JKU Business School, Johannes Kepler University Linz, Altenberger Straße 69, 4040 Linz, Austria}
\address[2]{Institute for Transport and Logistics Management, Vienna University of Economics and Business, Welthandelsplatz 1, 1020 Vienna, Austria}


\begin{abstract}
Many real-world optimisation problems involve multiple objectives. When considered concurrently, they give rise to a set of optimal trade-off solutions, also known as efficient solutions. These solutions have the property that neither objective can be improved without deteriorating another objective.
Motivated by the success of matheuristics in the single-objective domain, we propose a linear programming-based matheuristic for tri-objective binary integer programming. To achieve a high-quality approximation of the optimal set of trade-off solutions, a lower bound set is first obtained using the vector linear programming solver Bensolve. Then, feasibility pump-based ideas in combination with path relinking are applied in novel ways 
so as to obtain a high quality upper bound set. 
Our matheuristic is compared to a recently-suggested algorithm that is, to the best of our knowledge, the only existing matheuristic method for tri-objective integer programming. In an extensive computational study, we show that our method generates a better approximation of the true Pareto front than the benchmark method on a large set of tri-objective benchmark instances.
Since the developed approach starts from a potentially fractional lower bound set, it may also be used as a primal heuristic in the context of linear relaxation-based multi-objective branch-and-bound algorithms. 
\end{abstract}


\end{frontmatter}

\section{Introduction}\label{sec::intro}
Many real-world optimisation problems involve several and often conflicting objectives related to, e.g., costs, environmental impact and service quality. In a multi-objective optimisation problem, no single feasible solution simultaneously optimises all considered objectives, giving rise to a set of optimal trade-off solutions known as \textit{efficient} (or \emph{Pareto optimal}) solutions. Each efficient solution is the pre-image of a non-dominated point, when plotted in the space of objective values (the \emph{criterion space}). The set of all non-dominated points is called \emph{Pareto front} (see, e.g., the book \citep{ehrgott2005multicriteria} for an introduction to multi-objective optimisation).
In this study, we focus on tri-objective binary integer programming problems. 

To tackle a generic multi-objective integer programming (MOIP) problem, several exact algorithms that find one solution per non-dominated point have recently been proposed. Some of these work in the criterion space, such as, e.g., \citet{boland2017quadrant,tamby2021enumeration}, while others work in the space of feasible solutions, the \emph{decision space}, as initially suggested by \citet{kiziltan1983algorithm}. The latter rely on branch-and-bound based ideas.
\citet{forget2020branch,forget2022warm}, e.g., have recently developed branch-and-bound algorithms for problems with more than two objectives, which rely on linear programming (LP) relaxation-based bound sets \citep{ehrgott2007bound}.
However, using the available generic exact methods, only comparably small instances can be solved to optimality within reasonable computation times. This limitation motivates the development of a generic matheuristic approach for multi-objective integer programming that can obtain a high-quality approximation of the true Pareto front.

In branch-and-bound algorithms for single-objective mixed integer programming (MIP) problems, primal heuristics are often used to provide high-quality feasible solutions at an early stage, thus enabling to tighten the bound, resulting in finding an optimal solution more quickly. Hence, they have become an essential ingredient of state-of-the-art solvers such as CPLEX and Gurobi. In particular, LP relaxation-based primal heuristics have shown to be of value 
(see e.g., \citet{fischetti2003local}, \citet{fischetti2005feasibility})
Motivated by their success, our matheuristic first obtains a set of extreme points of the LP relaxation 
referred to as a lower bound (LB) set, using the vector linear programming solver \emph{Bensolve} \citep{lohne2017vector}. 
Depending on the fractionality of the LB set, we then either employ 
a variant of the \textit{feasibility pump} (\emph{FP}) or a combination of \textit{FP} and \textit{path relinking} (\emph{PR}) to generate \emph{feasible} integer solutions. We note that both \emph{FP} \citep{fischetti2003local} and \emph{PR} \citep{glover1997tabu} are generic heuristic techniques which were initially developed for single-objective optimisation.
Our method can stand alone as a matheuristic, or it can be used as a primal heuristic within LP relaxation-based tri-objective branch-and-bound algorithms.

The contributions of this study are as follows:
\begin{itemize}
  \item We propose a generic LP relaxation-based matheuristic algorithm for tri-objective binary integer programming, combining \emph{Bensolve} with \textit{FP} and \emph{PR}.
\item We propose a \emph{PR} scheme that generates integer solutions starting from fractional solutions. 
\item We show that generating additional fractional solutions from available (potentially integer) solutions to obtain new high quality integer solutions with \emph{FP} based ideas has a positive impact.
\item We propose two new sets of benchmark instances for multi-objective optimisation: One set consists of tri-objective facility location problems and the other set is based on instances from the well-known MIPLIB benchmark collection for single-objective problems \citep{miplib}. The instances are made available online.
\item We conduct an extensive computational study showing the efficacy of our generic algorithm compared to the only existing benchmark method.
\end{itemize}
We developed a new method built on the matheuristic designed for the multi-objective knapsack problem presented at the International Conference on Operations Research and Enterprise Systems (ICORES) 2021~\citep{ICORES-paper}.

The remainder of the paper is structured as follows. Section \ref{sec::related work} reviews related work, while Section \ref{sec::preliminaries} provides
the basic concepts and background for MOIP. 
Our matheuristic is described in
Section \ref{sec::algo}, and computational results are presented in Section \ref{sec::comp}. Finally, we provide concluding remarks in Section \ref{sec::conclusion}.

\section{Related work}\label{sec::related work}
A matheuristic is a hybrid approach combining mathematical programming with metaheuristics \citep{boschetti2009matheuristics}. Although matheuristics have been applied successfully in the single-objective domain (see e.g. \cite{Speranza_Archetti_2014}), comparably few contributions exist for MOIP problems. In the following, we discuss LP relaxation-based matheuristics in both the single and the multi-objective domains. Further, several successful studies using 
$PR$, the heuristic on which our method relies, in the multi-objective domain are presented.
\subsection{Matheuristics for single-objective optimisation}
Local branching, a combination of local search and MIP is firstly described by \citet{fischetti2003local} for single-objective MIP. Inspired by local search, they develop a method that can effectively explore the neighbourhood of solutions defined by linear cuts named local branching constraints. 
The proposed algorithm finds a new solution while solving a sub-problem made by fixing some variable values, and a k-opt operator.
For the given current feasible solution, its neighbourhood is built by performing soft fixing that maintains a certain number of decision variables in the incumbent solution but does not specifically fix any of those variables.
Secondly, a local branching constraint is introduced, to find a neighbouring solution within a certain Hamming distance. 
\citet{danna2005exploring} extend the idea of local branching and propose a relaxation induced neighbourhood search for single-objective MIP. Relying on the observation that an incumbent solution and a current LP relaxation solution often have common variable values, the method builds a promising neighbourhood by fixing them. Further, the defined sub-problem of the remaining variables is solved. The sub-problem of a relaxation induced neighbourhood search contains fewer variables than that of local branching and is solved faster. 
$FP$ has been introduced by \citet{fischetti2005feasibility} for single-objective MIP. 
The key idea of $FP$ is to build a sequence of roundings that possibly converges to a feasible solution of the given MIP. The algorithm requires a pair of points: an LP feasible but possibly fractional solution and its rounding that may be infeasible . At every iteration, the method attempts to minimise the distance between two points. While closing the distance between them, the algorithm may find a feasible solution. To be specific, once the rounded solution is obtained from the LP relaxation solution, the closest solution to it is found by solving an LP. If the new LP solution is a feasible solution to the MIP, the search stops. Otherwise, the LP solution is used in the next iteration. 
Since the algorithm performs well at finding feasible solutions (see e.g., \citet{achterberg2007improving,fischetti2009feasibility,boland2014boosting}), we use the principle of $FP$ for our method.

\subsection{Matheuristics for multi-objective optimisation}
\citet{soylu2015heuristic} designs a matheuristic algorithm for bi-objective mixed binary integer linear programming based on a variable neighbourhood search \citep{hansen2001variable} and local branching \citep{fischetti2003local}. The proposed method finds the approximating segments of the Pareto front and merges them at the end of the search. 
Another matheuristic framework for bi-objective binary integer programming has been proposed by \citet{leitner2016ilp}. The authors exploit that efficient solutions in the same region often have common characteristics. Informed by this fact, their algorithm finds feasible solutions by fixing a large number of variables based on existing solutions, similar to the relaxation induced neighbourhood search for single-objective optimisation. They also propose a multi-objective extension of local branching.
\citet{pal2019feasibility} develop an algorithm for bi-objective integer programming combining several existing algorithms in the literature of both single and bi-objective optimisation.
Their two-stage approach starts with the weighted sum method \citep{aneja1979bicriteria}, which is an \emph{inner approximation} algorithm to obtain an LB set for a bi-objective LP. Then they use the $FP$ algorithm with the LB set to generate feasible solutions. To improve the solutions obtained by the $FP$, they integrate the $FP$-based heuristic with a local search. This method is extended in \citet{pal2019fpbh} for multi-objective MIP. The authors develop a different weighted sum method for higher-dimensional problems and use a variant of local branching \citep{requejo2017feasibility} to increase the number of solutions. 
Since their extended algorithm, to the best of our knowledge, is the only generic matheuristic approach that can deal with MOIP problems with more than two objectives, we use it as a benchmark and refer to it as a feasibility pump-based heuristic ($FPBH$).
\citet{gandibleuxprimal} propose a primal heuristic computing an upper bound set (for minimisation problems) that can be embedded into a branch-and-bound algorithm for multi-objective binary integer programming.
LP solutions are first obtained by employing the $\epsilon$-constraint method \citep{haimes1971bicriterion}. Thereafter, inspired by $FP$, rounding is applied to the LP solution to find a feasible integer solution that should exist in the restricted criterion space defined by two adjacent images and the projection of the LP solution. The outcome set is not necessarily a close approximation of the Pareto front as their goal is to find a good upper bound set in a competitive computation time.

\subsection{Path relinking}

\textit{PR}, initially introduced by \citet{glover1997tabu} for single-objective problems, is a way of exploring trajectories between elite solutions in a neighbourhood space, and while exploring, the algorithm may find new solutions. Based on the fundamental idea that high-quality solutions share common characteristics, we expect a path between solutions might yield new solutions that share significant attributes inherited from the parent solutions. 
The procedure of \textit{PR} is as follows.
First an initiating and a guiding solution are chosen from a given set of solutions. These represent a starting point and destination, respectively. Starting from the initiating solution, the algorithm builds a certain number of neighbouring solutions in each step to reach the guiding solution. Among the created neighbouring solutions, one is selected as a new initiating solution towards which the algorithm moves. As the search proceeds, the new initiating solution receives more attributes of the guiding solution. The search ends when the initiating solution becomes identical to the guiding solution.
As $PR$ requires initial solutions, it can be intuitively applied to multi-objective optimisation problems where its LB set consists of a set of solutions. 

Starting from the study by \citet{gandibleux2003impact} who introduce $PR$ to multi-objective optimisation to solve bi-objective assignment problems, $PR$ has been used for diverse applications in combination with other metaheuristics.
For example, \citet{basseur2005path} use $PR$ with a genetic algorithm to solve bi-objective permutation flow-shop problems.
\citet{parragh2009heuristic} integrate $PR$ with a variable neighbourhood search to tackle bi-objective dial-a-ride problems. 
On the other hand, \citet{fernandes2021multi} employ a decomposition method which divides the criterion space into a set of sub-regions using a predefined set of direction vectors along with $PR$ to solve multi-objective combinatorial optimisation problems. 

\subsection{Our matheuristic}

Compared to the primal heuristic \citep{gandibleuxprimal} which aims to compute a good upper bound quickly, our matheuristic is designed to obtain a high-quality approximation of the true \textit{Pareto front}. Similar to~\cite{pal2019fpbh}, we initially construct an LB set by solving the LP relaxation of the MOIP. 
However, while \cite{pal2019fpbh} use a heuristic that generates a subset of the extreme points of the LB set by solving weighted sum problems, we use \emph{Bensolve} for that purpose, which relies on an \emph{outer approximation} algorithm~\cite{lohne2017vector}.  
Outer approximation algorithms generate valid LB sets even if stopped early. This can be useful in the context of multi-objective branch-and-bound algorithms. Thus our matheuristic could be more suitable as primal heuristic within multi-objective branch-and-bound algorithms compared to approaches relying on inner approximations methods like the weighted sum.
Then we apply a generalisation of the idea of the $FP$ framework \citep{pal2019feasibility} and a novel $PR$ scheme to the LB set (which can contain fractional solutions and integer solutions) in order to obtain (additional) feasible integer solutions.

\section{Preliminaries}
\label{sec::preliminaries}
Our matheuristic is designed to solve multi-objective binary integer programming (MOIP) problems.
In the following, we state a MOIP model in its general form, assuming minimisation of the objectives:

\begin{equation}
\label{eq:MO}
    \begin{aligned}
        Y = min\{Cx: Ax \geq b, x \in \{0,1\}^n \}, 
        \end{aligned}  \tag{MOIP}
\end{equation}\\
\noindent where $x_j$, $j=1,2,\dots,n$, is the vector of the decision variables. $X:= \{Ax \geq b, x \in \{0,1\}^n\}$ is the feasible set and $Y$ is the set of points in the criterion space, each of which corresponds to at least one solution vector $x \in X$. $C$ is a $p \times n$ objective function matrix where $c^k$, ($k=1,2,\dots,p$), is the $k^{th}$ row of $C$. In our case, $p=3$. $A$ is an $m \times n$ constraint matrix and $b$ is the right-hand-side vector for these constraints.

\subsection{Pareto dominance and supported solutions}
In multi-objective optimisation, a popular concept for determining the quality of a solution is $Pareto$ dominance \citep{ehrgott2005multicriteria}.
Suppose there are two solutions $x$ and $x^\prime$ to a problem (\ref{eq:MO}).
Then, the solution $x$ is called \textit{efficient/Pareto optimal} 
if and only if $c^k(x) \leq c^k(x')$ for all $k \in \{1,\dots,p\}$ and $c^k(x) < c^k(x')$ for at least one $k$. We state that $c^k(x)$ \textit{dominates} $c^k(x')$. Furthermore, 
$x$ is \textit{weakly efficient} 
if and only if $c^k(x) \leq c^k(x')$ for all $k\in\{1,\dots,p\}$ and we state that $c^k(x)$ \textit{weakly dominates} $c^k(x')$.
When a feasible solution $x$ is (\textit{weakly}) \textit{efficient}, $c^k(x)$ is called (\textit{weakly}) \textit{non-dominated}.
The \textit{efficient} set $X_E \subseteq X$ is defined as $X_E:=\{x \in X : \nexists \,x^\prime \in X : c(x^\prime) \leq c(x)\}$, and its image is referred to as the \textit{non-dominated} set $Y_N:=\{ c(x)\; | \;x \in X_E \}$.
The set of all \textit{non-dominated} points is called the \textit{Pareto front}.

If a solution $x$ can be found using the weighted sum method (i.e. optimising a convex combination of all the objective functions), it is called a \textit{supported efficient} solution. Otherwise, it is a \textit{non-supported efficient} solution.

\subsection{Lower bound set and LP relaxation}
The notion of a bound set was introduced by \citet{ehrgott2007bound}.
According to the definition established by the authors, an LB set is proposed to bound a subset $Y^\prime \subset Y$ of feasible points.
An LB set for $Y'$ is a subset $L \subseteq \mathbb{R}^{p}$ such that
\renewcommand{\theenumi}{\roman{enumi}$)$}%
\begin{enumerate}
  \item for each $y \in Y^\prime$ there is some $l\in L$ such that $l\leq y$, and
  \item there is no pair $y\in Y^\prime$, $l\in L$ such that $y$ \textit{dominates} $l$.
\end{enumerate}

One common way to obtain an LB set for a minimisation problem is to solve the LP relaxation of the original problem. 
In this study, for instance, binary variables become real variables bounded by 0 and 1 in the LP relaxation problem. We use this method as it allows us to compute the efficient set of the LP relaxation quickly.
In our case, we use the vector linear programming solver \textit{Bensolve} to generate the LB set.


\section{LP relaxation-based matheuristic}\label{sec::algo}
We propose an LP relaxation-based matheuristic which relies on \emph{Bensolve} to obtain the extreme points of the LB set.  In the following we call them \emph{LB solutions}.
Then, depending on the \emph{fractionality} of the LB set, we choose a different heuristic to generate feasible integer solutions. The fractionality represents the ratio between the number of fractional solutions and the total number of solutions in the LB set.
The first heuristic applies a variant of $FP$ while the second approach relies on a combination of $FP$ and $PR$. 
Both $FP$ and $PR$ take input solutions which are, in the general case, fractional. 
For problems which have a totally unimodular constraint matrix (e.g. the multi-objective assignment problem), the LB set only contains feasible integer solutions, i.e. it contains all supported efficient solutions. Moreover, it can also happen that by chance the fractionality of the LB set of an instance of a general multi-objective problem is quite low. In this case, the basic $FP$ cannot be used (efficiently) as it requires fractional solutions to generate feasible integer solutions. Thus, we develop a $FP$ variant, referred to as $FP^+$, that creates additional fractional solutions based on the current LB solutions from which the algorithm can start. These additional solutions are then used in $FP^+$ to generate new feasible integer solutions. 
%
%
For instances where most LB solutions are fractional, the basic $FP$ is applied first to find feasible integer solutions, then these are provided to $PR$ together with those LB solutions to which $FP$ was not applied (e.g. because of reaching a time limit). 

\subsection{General outline of our matheuristic}

We first give a general description of our matheuristic before we give details of the individual building blocks in the following subsections. The entire algorithmic framework is outlined in Algorithm \ref{algo::frame}.
In the first step of the proposed algorithm, LB solutions are obtained by running \emph{Bensolve} for at most the given time limit, \emph{BensolveTimeLimit} (line 1).
Depending on the fractionality of the LB solutions, we use two different approaches. If less than a certain percentage of the LB solutions, referred to as $allowedFractionality$, are fractional, we employ $FP^+$ (see Section \ref{sub:FP+}) that firstly produces additional fractional solutions then uses them to generate integer solutions, $X$ (lines 4-5).  
Otherwise, we use the basic $FP$ idea at first to provide feasible integer solutions together with the unused LB solutions for the $PR$ algorithm (lines 6-7). We refer to this method as \textit{feasibility pump generic path relinking} ($FPGPR$, see Section \ref{sub:FPGPR}).
Once $FPGPR$ receives the feasible integer solutions found by $FP$ and the unused LB solutions, the algorithm runs iteratively until it reaches its given time limit and returns $X$.
\begin{algorithm}
 \KwInput{LP relaxation problem of MOIP ($LP$)} 
 \KwOutput{approximate Pareto front $X$}
  Run \emph{Bensolve} for at most \emph{BensolveTimeLimit}\\ 
  $\Tilde{X} \leftarrow$ \emph{Bensolve}($LP$)\\
  Compute $fractionality$ of $\tilde{X}$\\
  \eIf{fractionality of $\Tilde{X} \leq $ allowedFractionality }
  {
  $X \leftarrow FP^+(\Tilde{X})$\\
  }
  {
    $X \leftarrow$\textit{FPGPR}($\Tilde{X}$)\\
  }
 \Return $X$
 \caption{LP relaxation-based matheuristic (\textit{LPBM})}
 \label{algo::frame}
\end{algorithm}


\subsection{Core$FP$}\label{sub:FP}
The $FP$ algorithm we use is based on the framework suggested by \citet{pal2019feasibility}, who extend the $FP$ idea to deal with bi-objective integer programming. 
In their algorithm, a maximum number of iterations is determined by the number of fractional values in the current LB solution. However, after initial experiments we fix 
the number of iterations per candidate solution to a small number,
this way many LB solutions can be used within a limited time budget. Another difference from their method is that we do not consider the dominance relationship between the newly found solution and the solutions stored in the archive during the search. We accept the new solution if it is feasible and does not exist in the incumbent solution set while \citet{pal2019feasibility,pal2019fpbh} only accept \textit{non-dominated} solutions.
We refer to this algorithm as Core$FP$, explained in Algorithm \ref{algo::FPprocess}, and use them for both $FP^+$ and $FPGPR$.
Both approaches take an LB solution $\tilde{x}$ as an input and maintain the $Tabu$ list to retain the infeasible  solutions found during the search. Differences include the ways of iterating the algorithm within the time limit and archiving new solutions. Details will follow in the next sections.

The first step of the algorithm is rounding the LP solution $\tilde{x}$ to make it integer (line 1). If the rounded solution $x_r$ is feasible and does not exist in the \textit{Archive} it is inserted into it and the search is done for the current $\tilde{x}$ (line 3).
In $FPGPR$, $\tilde{x}$ is removed from $\tilde{X}$ at this point (to return only unused LB solutions at the end). Otherwise, we check whether $\tilde{x}$ is a solution found previously (lines 4–5). If it is, the flip operation begins to generate another integer solution, which is referred to as $x_{flip}$ (line 6) (see \ref{apdix:flip} for the details of the operation). When no solution is found by the flip operation, a new iteration begins with the next candidate solution (lines 7–8). Otherwise, we examine $x_{flip}$'s feasibility and whether it does not exist in the current solution set. If it meets these conditions, it is added to \textit{Archive} and a new iteration begins (it is eliminated from $\tilde{X}$ in $FPGPR$) (lines 10–11).
If $x_r$ does not exist in $Tabu$, it is added to the list and we find a new LP solution that is closest to $x_r$ by solving the same optimisation problem as used by \citet{pal2019feasibility} (see \ref{apdix:lp}). The algorithm starts the new iteration with the newly found $\tilde{x}$ (lines 12–14). If no feasible solution is found by solving the LP, the search for the current $\tilde{x}$ terminates and the new iteration starts with the next $\tilde{x}$ (lines 15–16).
\begin{algorithm}[htbp]
\KwInput{$\Tilde{x},Archive,Tabu$}
\KwOutput{$Archive$}
    $x_r \gets$  Round($\tilde{x}$)\\    
    \eIf{$x_r$ is feasible \& ${x_r} \notin $ Archive}{
    $Archive = Archive \cup \{x_r\}$; \,\,\,break\\
    }{
    \eIf{$x_r \in Tabu$}
        {$x_{flip} \leftarrow Flip(Tabu,\tilde{x},x_r)$\\
        \eIf{ $x_{flip}= \emptyset$}
        { break
        }{
        \If{$x_{flip}$ is feasible \& $x_{flip} \notin $ Archive}
        { $Archive = Archive \cup \{x_{flip}\}$;\,\,\, break}
            }
        }
        { $Tabu = Tabu \cup \{x_r\}$\\
        $\tilde{x} \leftarrow FindnewLP(x_r)$; \\
        \eIf{$\tilde{x} = \emptyset$}{break}
        {go to line 1}
        }
    }
    \Return $Archive$
 \caption{$CoreFP$}\label{algo::FPprocess}
\end{algorithm}

\subsection{Feasibility pump$^+$ ($FP^+$)}\label{sub:FP+}
Algorithm \ref{algo:FFP} shows the $FP^+$ procedure designed for problems where the fractionality of computed LB solutions is below the given parameter $allowedFractionality$.
It takes as input the LB solutions $\hat{X}$ obtained from \emph{Bensolve}. These initial LB solutions are copied to the $Tabu$ list to prevent duplicate iterations in case the newly found solution already exist in $\Tilde{X}$ (line 3). 
To provide a fractional solution for the $FP$ process, we take the mean of two random solutions selected from the LB solutions. Thereby the new fractional solution inherits common variables from its parents and receives fractional variable values (line 6). The algorithm continues unless all sets of solution pairs are used (lines 7-10) within the given time limit, \emph{AlgoTimeLimit}. The new integer solution generated by \textit{CoreFP} is stored in $\Tilde{X}$, therefore, it can be used during the search. After the search terminates, $FP^+$ returns the filtered integer solution set $X$ (line 11).
\begin{algorithm}
\KwInput{$\Tilde{X},FP\_iter,AlgoTimeLimit$}
\KwOutput{$X$}
 $SolPair$: set for storing used $x_1$-$x_2$ pair\\
 $Tabu$: list for storing infeasible  solutions found during the search \\
 $Tabu \leftarrow \Tilde{X}$; \,\,\,\, $SolPair= \emptyset$; \,\,\,\,  \\
  \While{time $\leq$ AlgoTimeLimit}
  {$iter=0$; \\
  Choose $x_1,x_2$ from $\Tilde{X}$ at random; \,\,\,\,
  $\tilde{x} =  0.5\times x_1+ 0.5\times x_2$\\
    \While{ $\{x_1,x_2\} \notin SolPair$ \& iter $ \leq FP\_iter$ } 
   {    $\Tilde{X} \gets$ \textit{CoreFP}($\tilde{x},\Tilde{X},Tabu$)\\
    $iter=iter+1$
    }
    $SolPair = SolPair \cup \{\{x_1,x_2\}\}$\\
  }
    $X \leftarrow$  \textit{FilterIPsol}($\Tilde{X}$)

\Return $X$
\caption{ Feasibility pump$^+$ ($FP^+$) }\label{algo:FFP}
\end{algorithm}
We describe the $FP^+$ procedure for a tri-objective assignment problem with three tasks and agents (see the description of the benchmark problem in \ref{appendix}).
To keep  the example simple, we only consider the feasibility of solutions ignoring the dominance relationship. The list $\Tilde{X}$ provides the LB solutions of the problem.
\begin{equation*}
\Tilde{X}=
  \begin{bmatrix}
    [1& 0& 0& 0& 1& 0& 0& 0& 1]\\
    [0& 1& 0& 1& 0& 0& 0& 0& 1] \\
    [0& 0& 1& 0& 1& 0& 1& 0& 0] 
  \end{bmatrix}
  \label{LBset}
\end{equation*}
In Table \ref{table:FP}, the first column represents the set of solutions used for the iterations. The third column ($\tilde{x}$) shows the new fractional solutions obtained by combining the two selected solutions. The last column ($Tabu$) stores the discovered infeasible  integer solutions.
The asterisks mark a feasible integer solution.

Suppose $\Tilde{X[1]}$ and $\Tilde{X[2]}$ are chosen in the first repetition to produce a fractional solution $\tilde{x}$. The mean value of them is, $0.5 \times \Tilde{X[1]}+ 0.5\times \Tilde{X[2]}=$ (0.5 0.5 0 0.5 0.5 0 0 0 1). We assume its rounded solution $x_r$ is infeasible , thus $x_r$ is stored in $Tabu$ and the LP is solved to find the closest $\tilde{x}$. Since the operation cannot find a feasible LP solution and returns the empty set, we start a new iteration with a new pair of solutions.
From the two solutions \{$\Tilde{X[2]},\Tilde{X[3]}$\} we obtain $\tilde{x}=$ (0 0.5 0.5 0.5 0.5 0 0.5 0 0.5). Assuming its rounded value is infeasible , we solve the LP and obtain the new $\tilde{x}=$ (0 0 1 0 1 0 1 0 0) that already exits in $Tabu$. Therefore, the flip operation is conducted and we find the feasible solution (1 0 0 0 0 1 0 1 0). 
After adding this solution to $\Tilde{X}$ the search continues until it meets the stopping conditions. 

\begin{table}[htbp]
\centering
\begin{tabular}{|c|c|c|c|c|c|}
  \hline
     Set &iter & $\tilde{x}$ & $x_r$ & $x_{flip}$ & $Tabu$ \\
  \hline
   \{1,2\} & 1  & 0.5 0.5 0 0.5 0.5 0 0 0 1 & 1 1 0 1 1 0 0 0 1  & & 1 1 0 1 1 0 0 0 1  \\
  \cline{2-6}
  &2 & $\emptyset$ & & & 1 1 0 1 1 0 0 0 1 \\
  \hline 
  \{2,3\} & 1 & 0 0.5 0.5 0.5 0.5 0 0.5 0 0.5 & 0 1 1 1 1 0 1 0 1& &  1 1 0 1 1 0 0 0 1  \\ 
  \cline{2-6} 
  & 2 & 0 0 1 0 1 0 1 0 0 & 0 0 1 0 1 0 1 0 0 &1 0 0 0 0 1 0 1 0* &1 1 0 1 1 0 0 0 1\\
  \hline 
\end{tabular} 
\caption{\textit{Fractional feasbility pump} procedure}\label{table:FP} 
\end{table}



\subsection{Feasibility pump generic path relinking (FPGPR)}\label{sub:FPGPR}
Regarding problems for which more than $allowedFractionality$ of the LB solutions are fractional solutions, we combine the generic framework of \textit{PR} with $FP$ and name it $FPGPR$, which is described in Algorithm \ref{algo::FPGPR}.
For both the $FP$-part and the $PR$-part, we assign half of the time limit each.
At the first stage, the $FP$ process continues until there are no more solutions in the LB solutions or it reaches the given time limit (lines 2-7). The LP solution ($\tilde{x}$) is randomly chosen from $\Tilde{X}$ and used for \textit{Core FP} (line 6). All newly found integer solutions are stored in the list $X_{FP}$ and passed on to $GPR$ in addition to the remaining 
LB solutions $X_u$, once the search is completed.

$GPR$ requires an initial solution set $candX$, which consists of the generated integer solutions and unused LB solutions in the first stage of the algorithm 
(line 8).
The algorithm maintains two lists, namely $IGPair$ and $TabuS_i$, to keep track of the used pairs of initiating ($S_i$) and guiding ($S_g$) solutions and newly found solutions during the search. In the initial steps, we randomly choose $S_i$ and $S_g$ from $candX$ (line 11). Once $S_i$ and $S_g$ are determined, we check whether the pair of solutions has been used. If not, $GPR$ runs until $S_i$ becomes identical to $S_g$ (lines 12–25).
If a fractional solution is used as a guiding solution, the algorithm cannot terminate the search. Accordingly, we design the algorithm to start a new iteration when it
finds a 
neighbouring solution that has already been generated.
To generate the neighbourhood, we compare two solutions, $S_i$ and $S_g$, and find out those indices 
where the decision variables differ (line 13). This information is stored in the list $dif$ to be used in the $CreateNeighbours$ operation. 
 How neighbouring solutions are created is described in Subsection \ref{sec::createNB} in detail.
If no neighbouring solution is built, the algorithm stops and starts another iteration (lines 15–16). Otherwise, it checks the feasibility of the neighbouring solutions
and archives feasible ones in $candX$ (lines 17-20). Among
these 
neighbouring solutions, we choose the next $S_i$ in the $NextS_i$ operation (line 21), which is explained in Subsection \ref{sec::nextsi}. When the output of $NextS_i$ is not included in $candX$, it is inserted into $TabuS_i$ to prevent the algorithm from taking the same solution in the neighbourhood space (lines 22-23). At the end of every iteration, the current $[S_i,S_g]$ pair is stored in the list $IGPair$ (line 25). After filtering the fractional and dominated solutions, we obtain the set of feasible integer solutions $X$ (line 26).


  

\begin{algorithm}[htbp]
\KwInput{$\Tilde{X},FP\_iter,GPR\_iter,AlgoTimeLimit$}
\KwOutput{$X$}
  $X_{FP} = \emptyset$; \,\,\,\, $Tabu = \emptyset$; \,\,\,\,

   \While{$\Tilde{X}\neq \emptyset$ \& time $\leq$ AlgoTimeLimit$\times 0.5 $}
  {  
  $iter=0$\\
  Select $\tilde{x}$ from $\Tilde{X}$ at random\\
        \While{iter $\leq $ FP\_iter }
        {
             $X_{FP} \gets$ \textit{CoreFP}($\tilde{x},X_{FP},Tabu$)\\
              $iter=iter+1$
        }
  }
  \vskip 2mm
  $candX = X_{FP} \cup \Tilde{X}_u$; \,\,\,\, $IGPair = \emptyset$; \,\,\,\, $TabuS_i=\emptyset$; \,\,\,\,  \vskip 1mm
  
\While{time $\leq$AlgoTimeLimit $\times 0.5$}{
 $iter=0$;\\
  Choose $S_i$ \& $S_g$ from $candX$ at random\\
  \While{$S_i \neq S_g$ \& [$S_i,S_g$] $\notin$ IGPair \& $iter\leq GPR\_iter$}{
    $dif  \leftarrow$  index set of different variable values;\\
  $neighbours\leftarrow$ \textbf{CreateNeighbours}( $S_i,dif,TabuS_i$)\\
  \eIf{$|neighbours|=\emptyset$}{ break; }
  {
  
    \For{j=1,...,$|neighbours|$}{
    \If{neighbours[j] is feasible \& not in $ candX$ }
    {$candX$ = $candX \cup$ $\{neighbours$[$j$]\}}
    }
    $S_i \leftarrow$ \textbf{Next$S_i$}($neighbours$, $S_i$)\\
    \If{$S_i \notin candX$}{$TabuS_i$ = $TabuS_i \cup \{S_i\}$}
  }
  $iter$ = $iter+1$
  }
    $IGPair$ = $IGPair\cup \{[S_i,S_g]\}$;\\
 }
  $X \leftarrow$ \textit{Postprocessing}($candX$);\\
  \Return $X$
 \caption{Feasibility pump generic path relinking (\textit{FPGPR})}
 \label{algo::FPGPR}
\end{algorithm}

\subsubsection{CreateNeighbours\label{sec::createNB}}
The \textit{CreateNeighbours} operation is introduced to generate neighbouring solutions of the current $S_i$. Algorithm \ref{algo::CreateNB} shows how we construct the neighbourhood.
First, the number of neighbouring solutions to be generated is decided. If all the values in $S_i$ are integers, then the number of neighbouring solutions is the same as the number of different decision variable values between $S_i$ and $S_g$. For each index in $dif$ in which the indices of the differing variable values are stored, the corresponding decision variable value of $S_i$ is switched at a time. Specifically, if the current value is 1, it changes to 0, and vice versa (lines 3–7). 
When the variable value of the respective index of $S_i$ is fractional, the value is changed to 0 and 1 such that we remove the fractional value and eventually generate integer solutions. Both generated neighbouring solutions are archived in $neighbours$ (lines 8–10). Among the neighbourhood, we remove solutions explored in previous iterations stored in $TabuS_i$ (line 11).

\begin{algorithm}[htbp]
 \KwInput{$dif,TabuS_i$}
 \KwOutput{$neighbours$}
  $neighbours=\emptyset$ \\
  \For{i in dif}
  {
      \eIf{$S_i[i]$ is integer}
      {
            \eIf{$S_i[i]=1$}{$S_i[i] \leftarrow 0$}{$S_i[i]\leftarrow 1$}
           $neighbours= neighbours\cup \{S_i\}$;
      }
      {
      $S_i[i]\leftarrow1; \,\,\, neighbours= neighbours\cup \{S_i\}$;\\ 
      $S_i[i] \leftarrow 0; \,\,\, neighbours= neighbours \cup \{S_i\}$; 
      }
  }
  $neighbours = neighbours - TabuS_i$\\
  \Return $neighbours$
 \caption{$CreateNeighbours$}
 \label{algo::CreateNB}
\end{algorithm}

\subsubsection{$NextS_i$\label{sec::nextsi}}
To select the next initiating solution $S_i$ in the neighbourhood, we compare the improvement in the objective values of the neighbours. For this purpose, two matrices, $ratio\_table$ and $rank\_table$, are introduced in Algorithm \ref{algo::next}. The $ratio\_table$ records the ratio of the objective values of each neighbour to the current objective value $objS_i$, for each objective (lines 4–6). Further, each column of the $ratio\_table$ is ranked and the rankings are given to the $rank\_table$ (lines 7–9). The rankings of each neighbour are summed and passed to the $nb\_degree$ list (line 11). The neighbour that has the highest degree is chosen as the next $S_i$ (line 12).

\begin{algorithm}[htbp!]
 \KwInput{$neighbours,S_i$}
 \KwOutput{$S_i$}
  \textit{ratio\_table}: $|neighbours| \times p$ matrix\\
  \textit{rank\_table}: $|neighbours| \times p$ matrix\\
  \textit{nb\_degree} : list of size $|neighbours|$
  \vskip 2mm
  
  \For{i=1,...,$|neighbours|$ }
  {
    \For{j=1,...,$p$}
        {
        \textit{ratio\_table}[$i$][$j$] = $\cfrac{neighbours[i][j]}{objS_i[j]}$ 
        }
    }
  \For{i=1,...,$|neighbours|$}
    {
    \For{j=1,...,$p$}
            {
            \textit{rank\_table}[$i$][$j$] = rank of $neighbours$[$i$][$j$] in $j^{th}$ column
            }   
  }
  \vskip 2mm
    \For{i=1,...,$|neighbours|$}
    {
    \textit{nb\_degree}[$i$] = sum(\textit{rank\_table} $i^{th}$ column)\\
    }
    $S_i \leftarrow$ highest-degree solution in $neighbours$\\
  \Return $S_i$
 \caption{$NextS_i$}
 \label{algo::next}
\end{algorithm}

As an illustrative example, we briefly present the $GPR$ procedure for a tri-objective knapsack problem (see the description of the benchmark problem in \ref{appendix}). In this example, we suppose that the initiating and guiding solutions are (0 1 1 0.5) and (0.2 1 0 1), respectively, and there is one \textit{non-dominated} solution in the neighbourhood for the sake of simplicity.
We have four items to be added into a knapsack whose capacity is 20. The corresponding profits and weights are given in the matrix $P$ and $W$, respectively. \begin{equation*}\label{eq::profits}
P = 
\begin{pmatrix}
3 & 5 & 7 & 2\\
6 & 1 & 8 & 10\\
4 & 7 & 1 & 4 
\end{pmatrix} \;\;\;\;\;\;\;\; 
W = 
\begin{pmatrix}
8 & 6& 7 & 4
\end{pmatrix}
\end{equation*}\\
To create the neighbourhood, we check the common decision variable values between the two solutions. In this example, the second item is the only in common, therefore $dif=$[1,3,4]. Further, by switching one value whose index is stored in $dif$ at a time, we can generate a new solution. If the variable to be switched is fractional, we create two neighbouring solutions by changing it to 1 and 0, as the decision variables in our problem are binary. Table \ref{table:PR} shows the outcome of the process. In the first row, we generate four solutions since there are three different decision variable values and one of them is fractional. According to the profit matrix $P$, the corresponding objective values of each solution are [16,20,14], [6,6,9], [12,9,8] and [14,19,12].
As the first solution exceeds the knapsack capacity, the fourth solution (0 1 1 1) that has the second highest objective values becomes the new initiating solution (marked with * in Table \ref{table:PR}).
The same procedure applies to the next iteration. Between the two neighbouring solutions, the second solution is selected as the first solution is infeasible. 
At the third iteration, one solution is generated, which automatically takes the initiating solution role. Since the only created solution at the fourth iteration (0 1 0 1) has been explored previously, the search stops, and the last initiating solution (1 1 0 1) is stored in the integer solution archive.

\begin{table}[htbp]
\centering
\begin{tabular}{|c|c|c|c|c|c|}
  \hline
  Initiating & Guiding & \multicolumn{4}{|c|}{Neighbourhood} \\
  \hline
  0 \textbf{1} 1 0.5 & 0.2 \textbf{1} 0 1 &1 1 1 0.5 & 0 1 0 0.5 & 0 1 1 0 & 0 1 1 1* \\
  \cline{1-6}
  0  \textbf{1} 1 \textbf{1}  & 0.2 \textbf{1} 0 \textbf{1} &  1 1 1 1 & 0 1 0 1*  &\multicolumn{1}{c}{}\\ 
  \cline{1-4}
  0 \textbf{1 0 1} &0.2 \textbf{1 0 1} & 1 1 0 1* &\multicolumn{2}{c}{}\\
  \cline{1-3} 
  1 \textbf{1 0 1} & 0.2 \textbf{1 0 1} & 0 1 0 1& \multicolumn{2}{c}{}\\
  \cline{1-3}
\end{tabular}
\caption{\textit{Generic Path relinking} procedure}\label{table:PR} 
\end{table}

\section{Computational study}\label{sec::comp}
We evaluate the performance of our matheuristic approach, referred to as $LPBM$, focusing on the comparison with $FPBH$. 
\emph{Bensolve} is provided by \citet{lohne2017vector} at \url{http://www.bensolve.org/}.
$LPBM$ is implemented in the programming language Julia. For the benchmark algorithm, we use the Julia implementation of $FPBH$ (with the default setting) which is publicly available at \url{https://github.com/aritrasep/FPBH.jl}.
Since \emph{Bensolve} employs the $GLPK$ solver
to solve the mathematical models, it is used in all the algorithms. The exact MOIP solver, proposed by \citet{kirlik2014new} ($KS$), is also employed in the experiment to obtain the true Pareto front for comparison purposes, i.e. the $KS$ results are for reference only. All the experiments are carried out on a Quad-core X5570 Xeon CPU @2.93GHz with 48GB RAM. Throughout the paper, a maximisation objective function is converted into a minimisation one by multiplying it by~$-1$. Therefore, every result reported is for a minimisation problem.

In the following, we first describe the performance measures and benchmark problem instances used to conduct the computational experiments. Second, the results of the experiments are presented at the end of the section.
\subsection{Performance indicators}
In heuristic MOIP, the result is an approximation of the true Pareto front. To compare the approximation sets, several quality indicators have been proposed in the literature.
We use two metrics to evaluate the quality of our approximation set and compare it to a reference set.

\subsubsection{Hypervolume indicator}\label{hv}
The hypervolume (HV) indicator, introduced by \citet{zitzler1999multiobjective}, measures the volume of the dominated space of all the solutions contained in the approximation set. To calculate the dominated space, a reference point must be used.  Usually, a reference point is the “worst possible” point in the criterion space. In this study, we use the point (2,2,2). The computed figures of HV are normalised as follows.
Let $Y^{k}_N$ be a set of the $k^{th}$ objective values of the true Pareto front and $y\in \mathbb{R}^{p}$ be an arbitrary \textit{non-dominated} point obtained from a heuristic algorithm. Then, the normalised values of the obtained point are:
\begin{equation*}
    \begin{aligned}
        \frac{y^k-min(Y^{k}_N) }{max(Y^{k}_N)-min(Y^{k}_N)} && k=1,\dots,p.
    \end{aligned} 
\end{equation*}
Higher HV values indicate a better approximation.
We use the publicly available HV computing program provided by \citet{fonseca2006improved} at \url{http://lopez-ibanez.eu/hypervolume#intro} to obtain the HV values.

\subsubsection{Multiplicative unary epsilon indicator}
The multiplicative unary epsilon indicator, introduced by \citet{zitzler2003performance}, represents the distance between an approximation set $A$ and a reference set $R$. The $\epsilon$ is defined as the minimum factor necessary to multiply set $R$ to make the transformed \textit{non-dominated} points associated with set $R$ \textit{weakly dominated} by those of set $A$. All the objective values are normalised to [1,2]. With this metric, an indicator value greater than or equal to 1, and a lower value indicates a better approximation set. 
To compute the unary epsilon values, we use the performance assessment package provided by \citet{pisa} at \url{https://sop.tik.ee.ethz.ch/pisa/?page=assessment.php}.

\subsection{Benchmark problem instances}

The first two sets of instances, namely, the  multi-objective assignment problem (\texttt{MOAP}) and multi-objective knapsack problem (\texttt{MOKP}), are the same as those on which $FPBH$  is tested (the respective mathematical models are provided in the \ref{appendix}). 
Both instance sets were generated by \citet{kirlik2014new} and are publicly available at \url{http://home.ku.edu.tr/~moolibrary/}. The other two sets of instances, namely \texttt{TOFLP} and \texttt{TOMIPLIB} are new and provided for download at \url{https://www.jku.at/fileadmin/gruppen/132/Test_Cases/LPBM_instances.zip}.
More details on each instance set are given below:
\begin{itemize} 

\item \texttt{MOAP}: The set consists of 10 sub-classes containing 10 instances (100 instances in total). It is divided into instance classes based on the 
    the number of tasks which varies from 5 to 50 in increments of 5
\item \texttt{MOKP}: The set consists of 10 sub-classes containing 10 instances (100 instances in total).
    It is is divided into instance classes based on the 
    the number of items which varies from 10 to 100 in increments of 10

\item \texttt{TOFLP}: These are newly generated tri-objective facility location problem (\texttt{TOFLP}) instances based on the bi-objective single source capacitated facility location problem introduced by \citet{gadegaard2019bi}. 
The total coverage of customer demand is introduced as the third objective function.
The set consists of 12 sub-classes that include 10 instances. It is divided into instance classes based on the number of facilities which varies from 5 to 60 in increments of 5. The number of customers is double that of facilities 

\item \texttt{TOMIPLIB}: These are newly generated tri-objective integer programming instances based on single-objective instances
from the MIP library (MIPLIB 2017, \url{https://miplib.zib.de/}) \citet{miplib}. Initially, 90 single-objective binary integer programming instances categorised as \textit{easy} (solvable within one hour by a single-objective solver in default settings) were selected from the library. For these instances, coefficients for the second and third objective functions are then generated by shuffling the coefficients of the original objective function. However, initial computations showed that the resulting tri-objective instance were very hard. For the set \texttt{TOMIPLIB} used in our computational study, we selected the five instances for which \emph{Bensolve} can find the complete set of LB solutions within 10 minutes.

\end{itemize}

\subsection{Experimental results}
Following initial experiments, in our computational experiments, we use \emph{BensolveTimeLimit} set to ten minutes. The limit \emph{AlgoTimeLimit} is set to two minutes for the instance sets \texttt{MOAP}, \texttt{MOKP}, and \texttt{TOFLP}, while it is set to one hour for the more difficult instances set \texttt{TOMIPLIB}. The parameters \emph{allowedFractionality}, $FP\_iter$ and $GPR\_iter$ are set to 20\%, 10 and 20, respectively.

We report the following results for each algorithm: frac (\%), the fractionality of the LB set in case it varies across the instance set, the number of solutions ($|Y|$), the CPU time (sec), the percentage of the HV indicator value of the reference set (HV (\%)), and the unary epsilon indicator value ($\epsilon$). 

For the \texttt{MOAP}, \texttt{MOKP} and \texttt{TOFLP} the reference set is provided by $KS$.
As the exact method $KS$ cannot solve the \texttt{TOMIPLIB} instances within a couple of hours, we take the combined solution sets of $FPBH$ and $LPBM$ as the reference set.
The figures in Tables (\ref{tb:ap}–\ref{tb:kp}) show the average results over 10 test instances.
All the results are averaged over five runs for each instance.

The results on the \texttt{MOAP} are provided in Table \ref{tb:ap}. \emph{Bensolve}, denoted by $BEN$, already finds more solutions than $FPBH$ throughout all sub-classes while requiring far less computation time. Specifically, from the sub-class with 15 tasks (n$\geq$15), the number of solutions of \emph{Bensolve} is more than twice that of FPBH. Moreover, this difference increases as the instance becomes larger. Both the quantity and the quality of the solutions of \emph{Bensolve} are better than those of $FPBH$. For instance, they reach more than 99\% of the maximum HV throughout all the problem sets. Furthermore, the HV value increases as the problem size rises. By contrast, the highest HV value of $FPBH$ is 95.56\% in the second smallest problem class (n=10), and it decreases as the problem size increases. All the unary epsilon indicator values of \emph{Bensolve} are better than those of $FPBH$. As the solutions produced by \emph{Bensolve} are all \textit{supported}, we generate \textit{non-supported} solutions using $LPBM$ to improve the solution quality, especially for smaller instances. In fact, a set of \textit{supported efficient} solutions computed by \emph{Bensolve} already provides a high-quality approximation of the Pareto front, showing that the \texttt{MOAP} may not be the best suitable benchmark for LP relaxation-based MOIP heuristics.

\begin{sidewaystable}
\centering
 \begin{tabular}{crrrrrrrrrrrrrrrrrr}
 \hline 
    \multicolumn{5}{c}{$|Y|$} && \multicolumn{4}{c}{CPUtime(sec)}  & &\multicolumn{3}{c}{HV(\%)}& & \multicolumn{3}{c}{$\epsilon$} &\\ \cline{2-5} \cline{7-10} \cline{12-14} \cline{16-18}
    n  &  $KS$* & $FPBH$ & $BEN$ & $LPBM\dagger$ & & $KS$& $FPBH$ & $BEN$ &$LPBM$ && $FPBH$ & $BEN$ & $LPBM$ && $FPBH$ & $BEN$ &$LPBM$ \\ 
    \hline
5		&	14.1	&	6.7	&	7.5	&	\textbf{13.0}	&	&	0.1	&	\textbf{0.2}	&	$\sim$0.0	&	2.0	&	&	95.41	&	99.11	&	\textbf{99.96}	&	&	1.18	&	1.14	&	\textbf{1.05}	\\
10		&	176.8	&	20.9	&	39.0	&	\textbf{43.2}	&	&	10.7	&	\textbf{0.3}	&	$\sim$0.0	&	2.6	&	&	95.56	&	99.36	&	\textbf{99.43}	&	&	1.13	&	1.06	&	\textbf{1.06}	\\
15		&	674.9	&	41.0	&	83.1	&	\textbf{85.2}	&	&	92.5	&	\textbf{0.9}	&	0.07	&	3.0	&	&	94.95	&	99.48	&	\textbf{99.48}	&	&	1.12	&	1.06	&	\textbf{1.05}	\\
20		&	1860.5	&	63.4	&	161.3	&	\textbf{162.6}	&	&	359.1	&	\textbf{2.5}	&	0.22	&	5.3	&	&	94.42	&	99.54	&	\textbf{99.54}	&	&	1.12	&	1.05	&	\textbf{1.05}	\\
25		&	3567.8	&	89.3	&	253.1	&	\textbf{254.8}	&	&	872.2	&	\textbf{6.0}	&	0.60	&	7.7	&	&	94.46	&	99.69	&	\textbf{99.69}	&	&	1.11	&	1.04	&	\textbf{1.04}	\\
30		&	6181.3	&	140.7	&	379.4	&	\textbf{380.9}	&	&	1859.7	&	\textbf{17.1}	&	1.16	&	18.4	&	&	94.14	&	99.73	&	\textbf{99.75}	&	&	1.11	&	1.04	&	\textbf{1.04}	\\
35		&	8972.3	&	161.6	&	501.4	&	\textbf{504.0}	&	&	3285.6	&	\textbf{32.6}	&	2.08	&	33.6	&	&	94.09	&	99.77	&	\textbf{99.77}	&	&	1.11	&	1.03	&	\textbf{1.03}	\\
40		&	14679.7	&	241.9	&	699.1	&	\textbf{701.8}	&	&	6426.0	&	\textbf{72.0}	&	3.82	&	73.7	&	&	94.14	&	99.76	&	\textbf{99.76}	&	&	1.11	&	1.03	&	\textbf{1.03}	\\
45		&	17702.2	&	240.5	&	838.0	&	\textbf{841.2}	&	&	9239.0	&	\textbf{103.3}	&	6.10	&	105.3	&	&	93.88	&	99.79	&	\textbf{99.79}	&	&	1.11	&	1.03	&	\textbf{1.03}	\\
50		&	24916.8	&	335.5	&	1034.8	&	\textbf{1036.9}	&	&	15814.8	&	TL	&	9.74	&	TL	&	&	93.92	&	99.78	&	\textbf{99.78}	&	&	1.11	&	1.02	&	\textbf{1.02}	\\
\hline   
\end{tabular}
\caption{comparing the algorithm performance on \texttt{MOAP} for $p$=3, * indicates the optimal values, best heuristic values are in bold. TL=120sec.
$\dagger FP^+$ is used as the LB solutions are integer for all instances (frac = 0\%).
} \label{tb:ap}
\end{sidewaystable}

\begin{figure}
\begin{minipage}{.5\textwidth}
    \captionsetup{justification=centering}
    \includegraphics[width=1.0\linewidth]{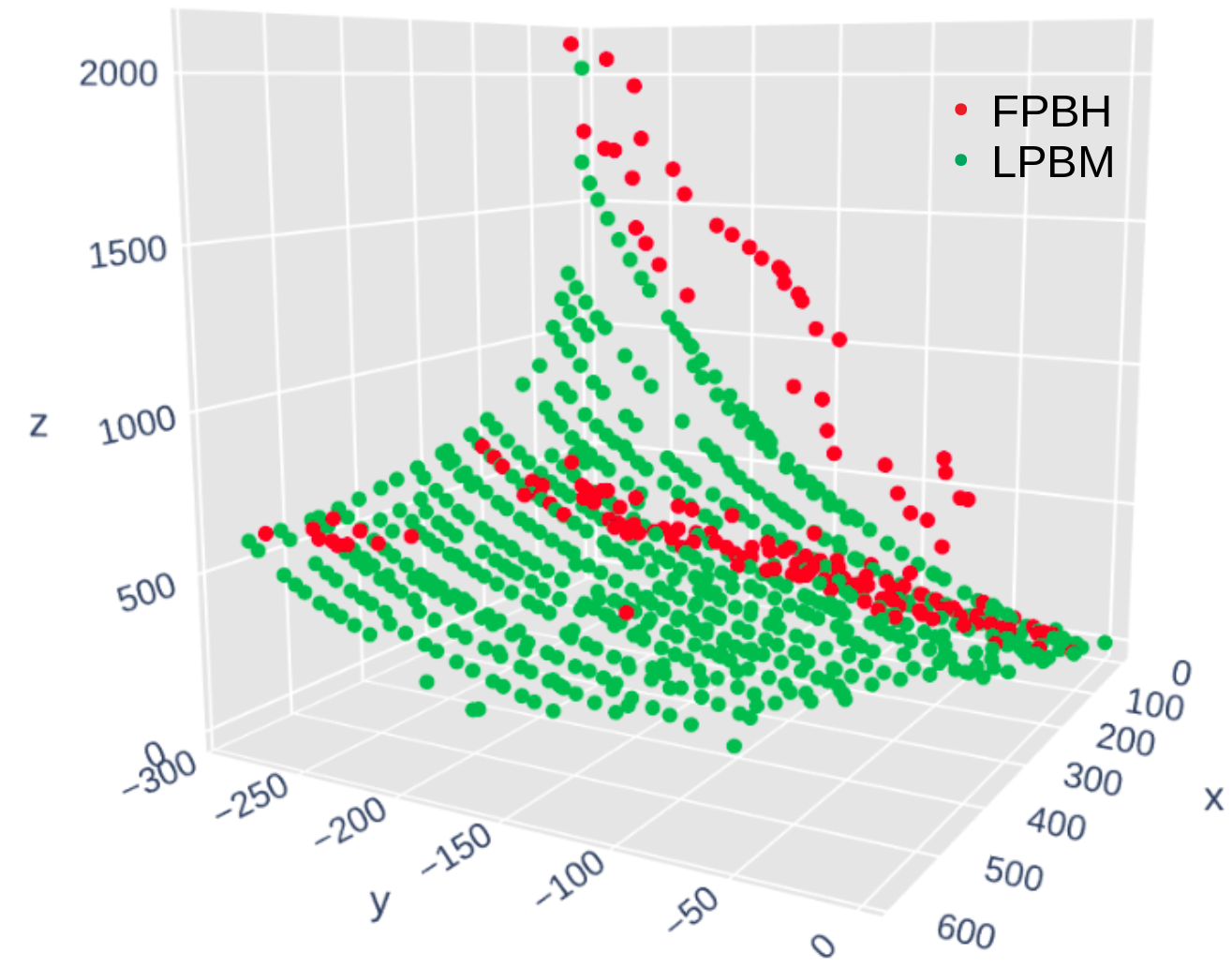}
    \caption{Generated solution points\\ for  TOFLP-n20\_40-ins3}
  \label{fig:FLPn2040}
\end{minipage}%
\begin{minipage}{.5\textwidth}
    \captionsetup{justification=centering}
  \includegraphics[width=1\linewidth]{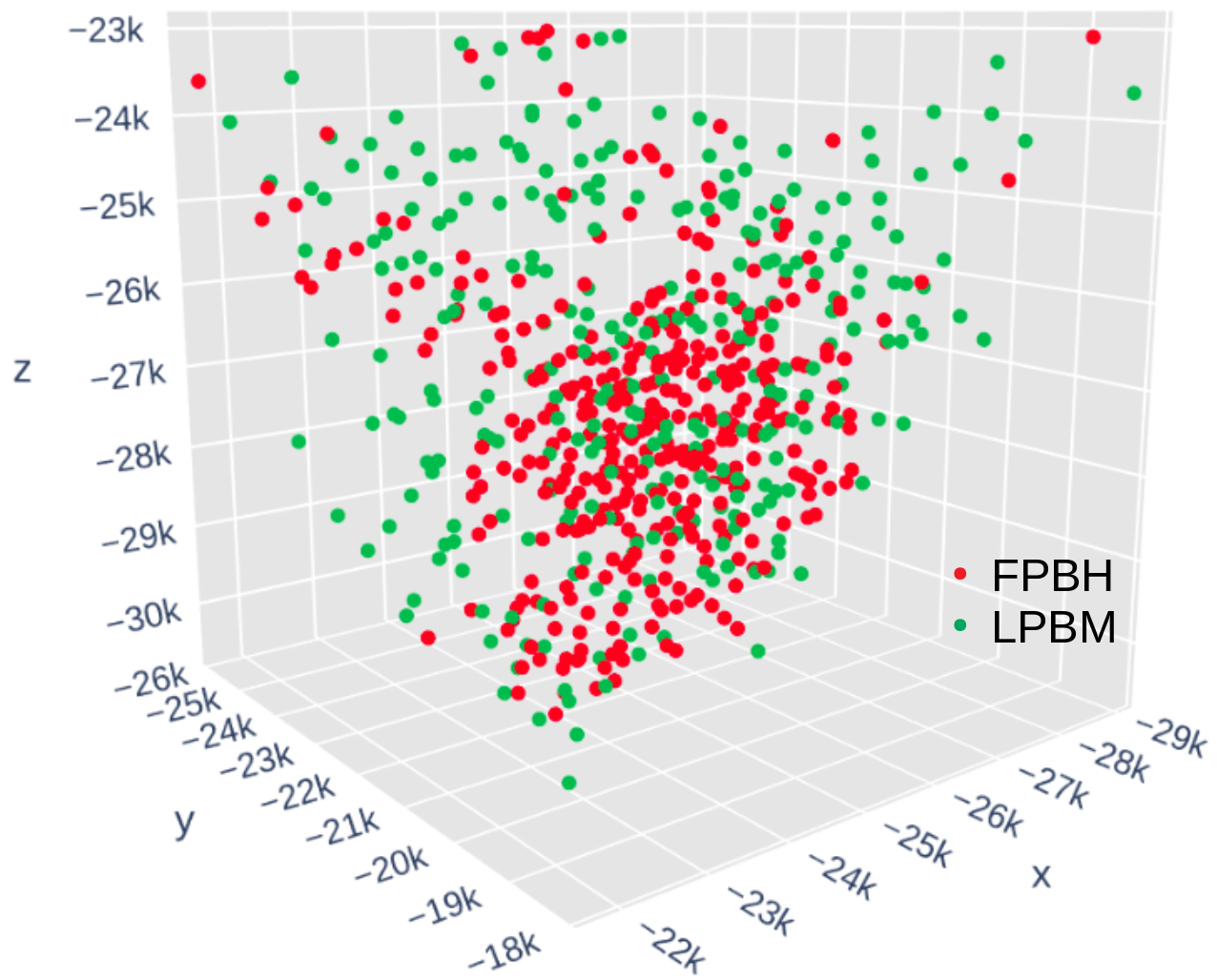}
  \caption{Generated solution points \\for MOKP-n70-ins3}
  \label{fig:KPn70}
\end{minipage}
\end{figure}

$LPBM$ results on the \texttt{MOKP} instances are provided in Table \ref{tb:kp}. Although $FPBH$ finishes computing over most of the \texttt{MOKP} instances and obtains more solutions for the larger sub-classes than $LPBM$, our method finds a higher-quality approximation for all the sub-classes based on both the HV and the unary epsilon indicators values. 
This can also be visually shown in Fig \ref{fig:KPn70}, which illustrates that the solution set of $LPBM$ is relatively well-distributed compared to that of $FPBH$.


\begin{sidewaystable}
\centering
 \begin{tabular}{crrrrrrrrrrrrrr}  
 \hline 
    &&\multicolumn{3}{c}{$|Y|$} & &\multicolumn{3}{c}{CPUtime(sec)}  & &\multicolumn{2}{c}{HV(\%)}& & \multicolumn{2}{c}{$\epsilon$} \\ \cline{3-5} \cline{7-9} \cline{11-12} \cline{14-15}
    n  & frac(\%) &  $KS$* & $FPBH$ & $LPBM\dagger$  & & $KS$& $FPBH$ & $LPBM$ & & $FPBH$ & $LPBM$ && $FPBH$ & $LPBM$ \\ 
    \hline
10	&	100.0&9.8	&	5.1	&	\textbf{7.6}	&	&	0.1	&	\textbf{0.2}	&	1.2	&	&	95.49	&	\textbf{97.52}	&	&	1.16	&	\textbf{1.12}	\\
20	&	100.0&38.0	&	18.2	&	\textbf{26.7}	&	&	1.0	&	\textbf{0.3}	&	1.4	&	&	96.20	&	\textbf{98.45}	&	&	1.12	&	\textbf{1.07}	\\
30	&	100.0&115.8	&	43.3	&	\textbf{55.9}	&	&	5.5	&	\textbf{0.6}	&	1.7	&	&	96.12	&	\textbf{98.21}	&	&	1.10	&	\textbf{1.06}\\
40	&	99.9&311.2	&	96.5	&	\textbf{102.2}	&	&	23.2	&	\textbf{1.5}	&	2.3	&	&	96.23	&	\textbf{98.54}	&	&	1.08	&	\textbf{1.05}\\
50	&	100.0&444.2	&	112	&	\textbf{117.5}	&	&	40.1	&	\textbf{2.7}	&	3.8	&	&	96.92	&	\textbf{98.48}	&	&	1.07	&	\textbf{1.05}\\
60	&	100.0&917.1	&	195.5	&	\textbf{200.8}	&	&	116.0	&	\textbf{6.4}	&	7.3	&	&	96.91	&	\textbf{98.65}	&	&	1.08	&	\textbf{1.04}	\\
70	&	99.9&1643.4	&	\textbf{346.6}	&	250.1	&	&	283.5	&	15.7	&	\textbf{15.3}	&	&	97.07	&	\textbf{98.62}	&	&	1.07	&	\textbf{1.04}	\\
80	&	100.0&2295.8	&	\textbf{441.5}	&	315.8	&	&	440.0	&	32.0	&	\textbf{31.3}	&	&	97.63	&	\textbf{98.80}	&	&	1.05	&	\textbf{1.04}	\\
90	&	99.9&3207.8	&	\textbf{503.6}	&	348.6	&	&	833.9	&	52.3	&	\textbf{50.3}	&	&	97.33	&	\textbf{98.79}	&	&	1.05	&	\textbf{1.03}	\\
100	&	99.9&5849.0	&	\textbf{894.2}	&	518.3	&	&	2478.4	&	\textbf{101.4}	&	119.1	&	&	97.18	&	\textbf{98.86}	&	&	1.05	&	\textbf{1.03}	\\

\hline   
\end{tabular}
\caption{comparing the algorithm performance on \texttt{MOKP} for $p$=3, * indicates the optimal values, best heuristic values are in bold. TL=120sec. $\dagger FPGPR$ is used as the fractionality of the LB set is above 20\% for all instances.} 
\label{tb:kp}
\end{sidewaystable}

Table \ref{tb:flp} shows  the results on the \textit{TOFLP} instances and a more noticeable performance difference between the algorithms. $FPBH$ terminates the search faster for small instances (n $\leq$15), but it finds far fewer solutions than $LPBM$ in all the sub-classes. In particular, for the sub-class with 30 items (n$=$30), the number of solutions generated by $FPBH$ starts to decrease and the gap between the two algorithms increases. A similar phenomenon is observed for instances with more than 25 items (n$\geq$25) for both quality measures; the performance of $FPBH$ declines, while that of $LPBM$ steadily improves as the instance size grows.
\begin{sidewaystable}
\centering
 \begin{tabular}{crrrrrrrrrrrrrr}  
 \hline 
      &&\multicolumn{3}{c}{$|Y|$} & &\multicolumn{3}{c}{CPUtime(sec)}  & &\multicolumn{2}{c}{HV(\%)}& & \multicolumn{2}{c}{$\epsilon$} \\ \cline{3-5} \cline{7-9} \cline{11-12} \cline{14-15}
    n  & frac(\%)  &$KS$* & $FPBH$ & $LPBM\dagger$  & & $KS$& $FPBH$ & $LPBM$ & & $FPBH$ & $LPBM$ && $FPBH$ & $LPBM$ \\ 
    \hline
5	&1.6	&120.3	&	19.8	&	\textbf{43.2}	&	&	3.4	&	\textbf{0.4}	&	2.4	&	&	97.64	&	\textbf{99.65}	&	&	1.13	&	\textbf{1.05}	\\
10	&0.5	&1360.9	&	59.0	&	\textbf{143.7}	&	&	180.0	&	\textbf{6.5}	&	9.7	&	&	97.86	&	\textbf{99.67}	&	&	1.17	&	\textbf{1.05}	\\
15	& 2.1&	3617.7	&	93.4	&	\textbf{329.3}	&	&	788.9	&	\textbf{58.0}	&	71.5	&	&	97.74	&	\textbf{99.69}	&	&	1.18	&	\textbf{1.04}	\\
20	&	4.0&8299.5	&	153.5	&	\textbf{520.6}	&	&	2752.0	&	TL	&	TL	&	&	98.26	&	\textbf{99.88}	&	&	1.18	&	\textbf{1.04}	\\
25	&	5.2&17502.8	&	179.0	&	\textbf{730.1}	&	&	8566.6	&	TL	&	TL	&	&	97.88	&	\textbf{99.90}	&	&	1.20	&	\textbf{1.03}	\\
30	&	4.9&25585.9	&	221.0	&	\textbf{966.7}	&	&	16891.7	&	TL	&	TL	&	&	96.39	&	\textbf{99.92}	&	&	1.23	&	\textbf{1.03}	\\
35	&	4.8&38791.5	&	159.0	&	\textbf{1176.9}	&	&	38859.4	&	TL	&	TL	&	&	95.87	&	\textbf{99.92}	&	&	1.24	&	\textbf{1.03}	\\
40	&	5.1&54791.4	&	115.6	&	\textbf{1420.7}	&	&	88433.7	&	TL	&	TL	&	&	95.84	&	\textbf{99.94}	&	&	1.24	&	\textbf{1.02}	\\
45	&	6.0&70778.5	&	80.5	&	\textbf{1668.8}	&	&	161139.3	&	TL	&	TL	&	&	94.89	&	\textbf{99.94}	&	&	1.24	&	\textbf{1.02}	\\
50	&	7.7&100187.1	&	64.2	&	\textbf{1943.8}	&	&	226236.8	&	TL	&	TL	&	&	95.56	&	\textbf{99.94}	&	&	1.24	&	\textbf{1.02}	\\
55	&	8.7&120342.7	&	58.8	&	\textbf{2152.0}	&	&	465427.9	&	TL	&	TL	&	&	95.67	&	\textbf{99.95}	&	&	1.23	&	\textbf{1.02}	\\
60	&	7.9 &148529.8	&	55.4	&	\textbf{2477.4}	&	&	804132.5	&	TL	&	TL	&	&	95.41	&	\textbf{99.95}	&	&	1.25	&	\textbf{1.02}	\\
\hline   
\end{tabular}
\caption{Comparing the algorithm performance on \texttt{TOFLP}, * indicates the optimal values, best heuristic values are in bold. TL=120sec. $\dagger FP^+$ is used as the fractionality of the LB set is below 20\% for all instances.} 
\label{tb:flp} 
\end{sidewaystable}

\begin{figure}[htbp!]
    \includegraphics[width=1.0\textwidth]{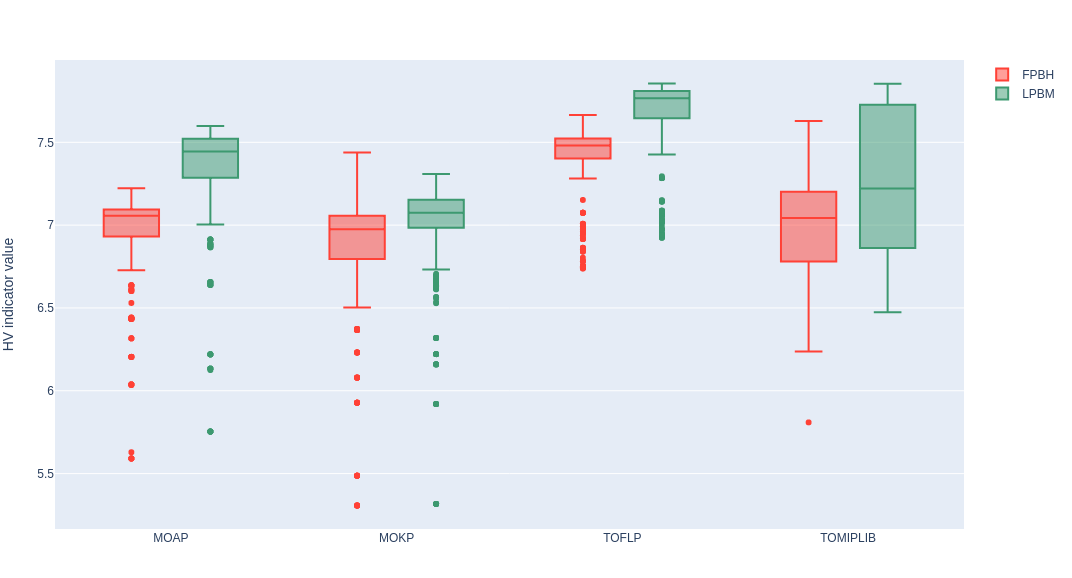}
  \caption{Comparing the HV values by instance set. (larger is better)}
  \label{fig:hv}
  \includegraphics[width=1.0\linewidth]{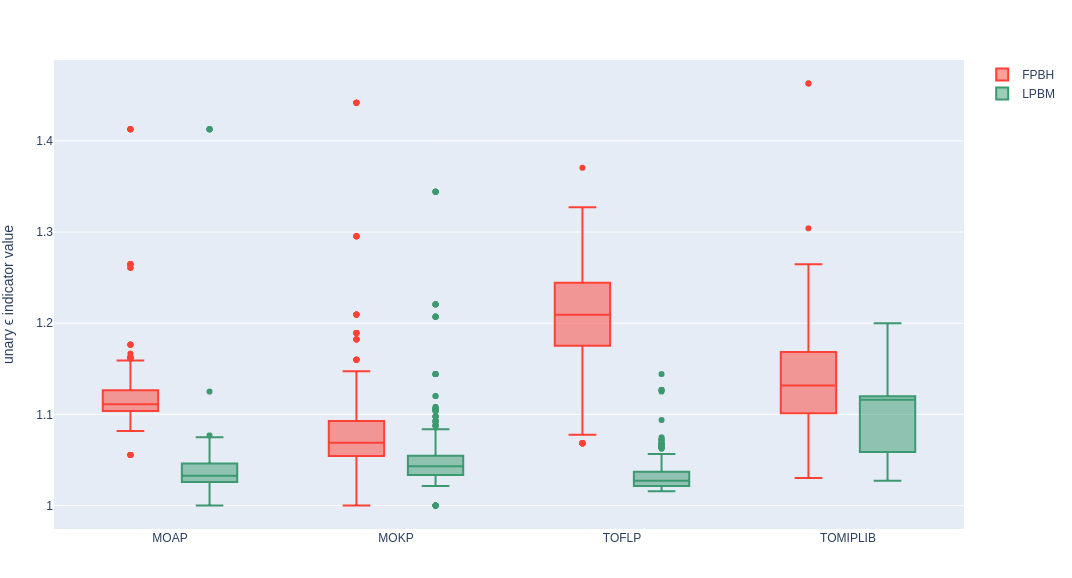}
  \caption{Comparing the unary $\epsilon$ values by instance set (smaller is better)}
  \label{fig:eps}
\end{figure}

Table \ref{tb:mip} shows the results on the \texttt{TOMIPLIB} instances. $FPBH$ and $LPBM$ are competitive in every aspect for set packing instances (cvs08r139,cvs16r70).
For the remaining instances, $FPBH$ finds considerably fewer solutions than $LPBM$. 
For \textit{neos} instances in which $FPBH$ takes less computation time, we conduct extra experiments on $LPBM$ by imposing the time that $FPBH$ took. In the experiments, we find that our method still outperforms $FPBH$. This results can be found in Table \ref{appendixtb}. The overview of the solution quality of $FPBH$ and $LPBM$ is given in Figures \ref{fig:hv} and \ref{fig:eps}. The comparison is made by calculating the indicator values on all the instances of each problem class.
\begin{table}[!ht]
\makebox[\linewidth]{
 \begin{tabular}{crrrrrrrrrrrr}  
 \hline 
    &\multicolumn{2}{c}{$|Y|$} & &\multicolumn{2}{c}{CPUtime(sec)}  & &\multicolumn{2}{c}{HV(\%)}& & \multicolumn{2}{c}{$\epsilon$} \\ 
    \cline{2-3} \cline{5-6} \cline{8-9} \cline{11-12}
    instance  &  $FPBH$ & $LPBM\dagger$  & &  $FPBH$ & $LPBM$ & & $FPBH$ & $LPBM$ && $FPBH$ & $LPBM$ \\
    \hline
cvs08r139	& 	\textbf{19.8}	&	\textbf{19.8}	&	&	\textbf{3477.9}	&	TL	&	& \textbf{94.75}	&	91.22	&	&	\textbf{1.17} &	1.18\\
cvs16r70	&	8.5	&	\textbf{9.0}	&	&	\textbf{3476.6}	&	TL	&	&95.37	&	\textbf{96.65}	&	&	1.15	&	\textbf{1.12}\\
n2seq36f	&	33.5	&	\textbf{132.8}	&	&	TL	&	TL	&	&	93.35	&	\textbf{99.56}	&	&	1.12	&	\textbf{1.05}\\
neos1516309	&	17.5	&	\textbf{306.2}	&	&	\textbf{1272.6}	&	TL	&	&87.79	&	\textbf{96.19}	&	&1.21	&	\textbf{1.11}\\
neos1599274	&	24.0	&	\textbf{108.8}	&	&	\textbf{1032.6}	&	TL	&	&94.78	&	\textbf{99.37}	&	&1.12	&	\textbf{1.07}\\
\hline   
\end{tabular}
}
\caption{comparing the algorithm performance on \texttt{TOMIPLIB}, best heuristic values are in bold. TL=3600sec. $\dagger FPGPR$ is used as the LB solutions are fractional for all instances (frac = 100\%).} 
\label{tb:mip}
\end{table}

\begin{table}[!ht]
\makebox[\linewidth]{
 \begin{tabular}{crrrrrrrrrrrr}  
 \hline 
    &\multicolumn{2}{c}{$|Y|$} & &\multicolumn{2}{c}{CPUtime(sec)}  & &\multicolumn{2}{c}{HV(\%)}& & \multicolumn{2}{c}{$\epsilon$} \\ 
    \cline{2-3} \cline{5-6} \cline{8-9} \cline{11-12}
    instance  &   $FPBH$ & $LPBM\dagger$  & &  $FPBH$ & $LPBM$ & & $FPBH$ & $LPBM$ && $FPBH$ & $LPBM$ \\
    \hline
neos1516309	&	17.5	&	\textbf{67.5}	&	&	1272.6 &\textbf{1220.7}	&	&91.38	&	\textbf{94.80}	&	&	0.148	&	\textbf{0.138}\\
neos1599274	&	24.0	&	\textbf{26.0}	&	&	\textbf{1032.6}	&	1035.3	&	&96.47	&	\textbf{98.07}	&	&	0.127	&	\textbf{0.111}\\
\hline   
\end{tabular}
}
\caption{Comparing the performance of the algorithms on the \texttt{TOMIPLIB}, best heuristic values are in bold. $\dagger FPGPR$ is used as the LB solutions are fractional for all instances (frac = 100\%).} 
\label{appendixtb}
\end{table}

\section{Conclusion}\label{sec::conclusion}
We develop as linear programming-based matheuristic for tri-objective binary integer programming. The proposed algorithm employs a high-performing vector linear programming solver, \emph{Bensolve}, to obtain a lower bound set. Further, depending on the fractionality of the lower bound set, different combinations of feasibility pump-based ideas and path relinking are used to generate integer solutions.
Since the proposed algorithm is designed to generate integer solutions starting from the lower bound set that might contain fractional values, it can be used to deal with any  multi-objective binary integer programming problem.
Our extensive computational study shows the efficacy of our algorithm which outperforms the benchmark method, $FPBH$, generating a high-quality approximation of the true \textit{Pareto front} in the vast majority of the instances.

For future work, we plan to tailor the algorithm to real-world applications such as in supply chain network design. 
Since the size of real-world problems is much larger, one major task could be reducing the computation time by employing a learning algorithm.

\section*{Acknowledgements}
This research was funded in whole, or in part, by the Austrian Science Fund (FWF) [P 31366-NBL] and [P 35160-N]. For the purpose of open access, the author has applied a CC BY public copyright licence to any Author Accepted Manuscript version arising from this submission.
\clearpage
\appendix

\section{The mathematical models of our benchmark problems}\label{appendix}\label{sec::problem}
\noindent The following problems are used to compare the performance of the algorithms. 

\subsection{Multi-objective assignment problem}
\noindent In the well-known assignment problem, a certain number of tasks $l \in \{1,\dots,n\}$ and agents $r \in \{1,\dots,n\}$ are given. The decision variable $x_{rl}$ takes the value of 1 if task $l$ is assigned to the agent and 0,  otherwise. The costs $c^1_{rl},\dots,c^p_{rl}$ are incurred if a task is allocated to an agent. The \texttt{MOAP} can be stated as follows:
\begin{align}
    \min \quad  \sum_{r=1}^{n}\sum_{l=1}^{n}{c_{rl}^{j}x_{rl}} &&& j=1,\ldots,p 
    \label{eq:AP1}
    \\
    \textrm{s.t.} \quad \sum_{l=1}^{n}{x_{rl}} = 1 &&& r = 1,\dots,n 
    \label{eq:AP(2)}
    \\
    \sum_{r=1}^{n}{x_{rl}}  = 1 &&& l = 1,\dots,n 
    \label{eq:AP(3)}
    \\
    x_{rl}  \in \{0,1\}  &&& r,l=1,\ldots,n.
    \label{eq:AP4}
\end{align}
\\
The goal of the \texttt{MOAP} is to find an optimal assignment of all the tasks to agents while minimising the $p$ cost functions (\ref{eq:AP1}). Equation (\ref{eq:AP(2)}) proposes a limit whereby each agent is assigned to only one task. Equation (\ref{eq:AP(3)}) ensures that each task is assigned to one agent only.

As the constraint matrix of \texttt{AP} is totally unimodular, every vertex of the LP relaxation is an integer vector. Thus, we can naturally obtain integer solutions of the \texttt{AP} by solving the LP relaxation. 

\subsection{Multi-objective knapsack problem}
\noindent In the \texttt{MOKP}, a set of items with a certain weight $w_r$ and profit $v_r$ is given. A decision maker must select a subset of these items such that the total weight does not exceed a given capacity $W$. The binary decision variable $x_r$ has a value of 1, when item $r$ is selected for the knapsack. Otherwise, $x_r$ is 0. Here, $W$, $v_r$, and $w_r$ are non-negative integer values. 
The \texttt{MOKP} model is stated as follows:

\begin{equation}
\label{eq:KP(1)}
\begin{aligned}
    &&\max \sum_{r=1}^{n}{v_{r}^{j}x_{r}} && j=1,\ldots,p \\
\end{aligned}  
\end{equation}
\begin{equation}
\label{eq:KP(2)}
\begin{aligned}
    &\textrm{s.t.} \sum_{r=1}^{n}{w_{r}x_{r}} \leq W \\
\end{aligned}  
\end{equation}
\begin{equation}
\label{eq:KP(3)}
\begin{aligned}
    &x_{r}\in \{0,1\}  & r=1,\ldots,n 
\end{aligned}  
\end{equation}\\
Equation (\ref{eq:KP(1)}) represents the objective functions maximising the $p$ total profits of selected items. Equation (\ref{eq:KP(2)}) denotes the capacity constraint. The total weight of the items stored in the knapsack cannot exceed the given capacity.

\subsection{Tri-objective facility location problem}
\noindent A set $I$ of potential facility sites and a set $J$ of customers are given. When a facility is opened, a fixed opening cost $f_i$ is incurred. A servicing cost $c_{ij}$ arises when the facility $i$ services a customer $j$ whose demand is represented by $d_j$. The bi-objective single source capacitated facility location problem \texttt{BO-SSCFLP} 
minimises the total servicing cost and fixed opening cost under the condition of each facility's capacity and customer demand. 
In the extended model, the \texttt{TOFLP}, we consider maximising the service coverage. Hence, the capacity limit of facilities is eliminated from \texttt{BO-SSCFLP} and the total coverage of customer demand is introduced as the third objective function to the \texttt{TOFLP}.
\begin{align}
\min (\sum_{i\in I}\sum_{j\in J} c_{ij} x_{ij}, \hspace{0.4cm} & \sum_{i \in I} f_i y_i, \hspace{0.4cm} -\sum_{j \in J} d_j  z_j )  \label{eq:FLP1} \\
\text{s.t.} \sum_{i \in I} x_{ij} & = z_j, 
&& \forall j \in J \label{eq:FLP2} \\
    x_{ij} & \leq y_i,     && \forall i \in I,j \in J  \label{eq:FLP3} \\
    x_{ij},y_i,z_j & \in \{0,1\}, 
    &&\forall i \in I, j \in J    \label{eq:FLP4} 
\end{align}

Equation (\ref{eq:FLP1}) defines three objectives: total servicing cost, total fixed opening cost, and total customer demand coverage. 
Equation (\ref{eq:FLP2}) associates serviced customers to demand coverage. 
Equation (\ref{eq:FLP3}) imposes that a customer should only be serviced by an opened facility.

\section{Feasibility pump components}
\subsection{Flip operation}\label{apdix:flip}
Algorithm \ref{algo:flip} describes the $Flip$ operation used in $LPBM$.
We refer to \citet{pal2019feasibility} for detailed explanations.

\begin{algorithm}[htbp]
 \KwInput{$Tabu,\tilde{x},x_r$}
 \KwOutput{$x_I$}
 N = dimension of $\tilde{x}$\\
 $\{e^1,\dots,e^N\} \leftarrow$ sort the set $j=\{1,\dots,N$\} based on the value of $|\tilde{x}^j-x_r^j|$ in descending order\\
 $x_I = Null; \,\,\,\, \hat{x}_r \leftarrow x_r; \,\,\,\, i \leftarrow 1$\\
 
  \While{$i\leq N$ and $x_I=Null$}
  {
    \eIf{$\hat{x}_r[e[i]]==1$}
      {$\hat{x}_r[e^i]=0$}
      {$\hat{x}_r[e^i]=1$}
    \eIf{$\hat{x}_r \notin Tabu$}
      {$x_I \leftarrow \hat{x}_r$}
      {$i=i+1$}
  }
  \If{$x_I==Null$}
    {$i\leftarrow 1$\\
    \While{$i\leq N$}
        {$\hat{x}_r \leftarrow x_r$\\
        $Num \leftarrow RandomBetween(\lceil N/2 \rceil, N-1)$\\
        $L \leftarrow RandomlyChoose(Num,\{1,\dots,N\})$\\
        \For{$l \in L$}
            {
            \eIf{$\hat{x}_r[e^l]==1$}
              {$\hat{x}_r[e^l]=0$}
              {$\hat{x}_r[e^l]=1$}
            }
        \eIf{$\hat{x}_r \notin Tabu$}
            {$x_I \leftarrow \hat{x}_r$}
            {$i=i+1$}    
        }
    }
  \Return $x_I$
 \caption{$Flip$}\label{algo:flip}
\end{algorithm}

\subsection{Optimisation problem for $FindnewLP$}\label{apdix:lp}
$FindnewLP$ updates $\tilde{x}$ using $x_r$ by solving the following optimisation problem,

\begin{equation*}
    \begin{aligned}
        \tilde{x} = \arg \min\left\{ \sum_{j=1:x_{r}^j=0}^{n}x^j + \sum_{j=1:x_{r}^j=1}^{n}(1-x^j): x \in \mathscr{X}
        \right\},
    \end{aligned}
\end{equation*}
where $n$ is the dimension of $x_r$ and $\mathscr{X}$ is the feasible set of LP relaxation solutions of MOIP.
We aim to find the LP solution that is closest to the current $x_r$ by solving the optimisation problem.


\clearpage
\bibliography{main}

\begin{thebibliography}{38}
\expandafter\ifx\csname natexlab\endcsname\relax\def\natexlab#1{#1}\fi
\providecommand{\url}[1]{\texttt{#1}}
\providecommand{\href}[2]{#2}
\providecommand{\path}[1]{#1}
\providecommand{\DOIprefix}{doi:}
\providecommand{\ArXivprefix}{arXiv:}
\providecommand{\URLprefix}{URL: }
\providecommand{\Pubmedprefix}{pmid:}
\providecommand{\doi}[1]{\href{http://dx.doi.org/#1}{\path{#1}}}
\providecommand{\Pubmed}[1]{\href{pmid:#1}{\path{#1}}}
\providecommand{\bibinfo}[2]{#2}
\ifx\xfnm\relax \def\xfnm[#1]{\unskip,\space#1}\fi
\bibitem[{Achterberg \& Berthold(2007)}]{achterberg2007improving}
\bibinfo{author}{Achterberg, T.}, \& \bibinfo{author}{Berthold, T.}
  (\bibinfo{year}{2007}).
\newblock \bibinfo{title}{Improving the feasibility pump}.
\newblock {\it \bibinfo{journal}{Discrete Optimization}\/},  {\it
  \bibinfo{volume}{4}\/}, \bibinfo{pages}{77--86}.
\bibitem[{An et~al.(2021)An, Parragh, Sinnl \& Tricoire}]{ICORES-paper}
\bibinfo{author}{An, D.}, \bibinfo{author}{Parragh, S.},
  \bibinfo{author}{Sinnl, M.}, \& \bibinfo{author}{Tricoire, F.}
  (\bibinfo{year}{2021}).
\newblock \bibinfo{title}{A {LP} relaxation based matheuristic for
  multi-objective integer programming}.
\newblock In {\it \bibinfo{booktitle}{Proceedings of the 10th International
  Conference on Operations Research and Enterprise Systems (ICORES 2021)}\/}
  (pp. \bibinfo{pages}{88--98}).
\bibitem[{Aneja \& Nair(1979)}]{aneja1979bicriteria}
\bibinfo{author}{Aneja, Y.~P.}, \& \bibinfo{author}{Nair, K.~P.}
  (\bibinfo{year}{1979}).
\newblock \bibinfo{title}{Bicriteria transportation problem}.
\newblock {\it \bibinfo{journal}{Management Science}\/},  {\it
  \bibinfo{volume}{25}\/}, \bibinfo{pages}{73--78}.
\bibitem[{Archetti \& Speranza(2014)}]{Speranza_Archetti_2014}
\bibinfo{author}{Archetti, C.}, \& \bibinfo{author}{Speranza, M.~G.}
  (\bibinfo{year}{2014}).
\newblock \bibinfo{title}{A survey on matheuristics for routing problems}.
\newblock {\it \bibinfo{journal}{EURO Journal on Computational
  Optimization}\/},  {\it \bibinfo{volume}{2}\/}, \bibinfo{pages}{223--246}.
\bibitem[{Basseur et~al.(2005)Basseur, Seynhaeve \& Talbi}]{basseur2005path}
\bibinfo{author}{Basseur, M.}, \bibinfo{author}{Seynhaeve, F.}, \&
  \bibinfo{author}{Talbi, E.-G.} (\bibinfo{year}{2005}).
\newblock \bibinfo{title}{Path relinking in pareto multi-objective genetic
  algorithms}.
\newblock In {\it \bibinfo{booktitle}{International Conference on Evolutionary
  Multi-Criterion Optimization}\/} (pp. \bibinfo{pages}{120--134}).
\newblock \bibinfo{organization}{Springer}.
\bibitem[{Bleuler et~al.(2003)Bleuler, Laumanns, Thiele \& Zitzler}]{pisa}
\bibinfo{author}{Bleuler, S.}, \bibinfo{author}{Laumanns, M.},
  \bibinfo{author}{Thiele, L.}, \& \bibinfo{author}{Zitzler, E.}
  (\bibinfo{year}{2003}).
\newblock \bibinfo{title}{{PISA} --- a platform and programming language
  independent interface for search algorithms}.
\newblock In \bibinfo{editor}{C.~M. Fonseca}, \bibinfo{editor}{P.~J. Fleming},
  \bibinfo{editor}{E.~Zitzler}, \bibinfo{editor}{K.~Deb}, \&
  \bibinfo{editor}{L.~Thiele} (Eds.), {\it \bibinfo{booktitle}{Evolutionary
  Multi-Criterion Optimization {(EMO 2003)}}\/} Lecture Notes in Computer
  Science (pp. \bibinfo{pages}{494--508}).
\newblock \bibinfo{address}{Berlin}: \bibinfo{publisher}{Springer}.
\bibitem[{Boland et~al.(2017)Boland, Charkhgard \&
  Savelsbergh}]{boland2017quadrant}
\bibinfo{author}{Boland, N.}, \bibinfo{author}{Charkhgard, H.}, \&
  \bibinfo{author}{Savelsbergh, M.} (\bibinfo{year}{2017}).
\newblock \bibinfo{title}{The quadrant shrinking method: A simple and efficient
  algorithm for solving tri-objective integer programs}.
\newblock {\it \bibinfo{journal}{European Journal of Operational Research}\/},
  {\it \bibinfo{volume}{260}\/}, \bibinfo{pages}{873--885}.
\bibitem[{Boland et~al.(2014)Boland, Eberhard, Engineer, Fischetti, Savelsbergh
  \& Tsoukalas}]{boland2014boosting}
\bibinfo{author}{Boland, N.~L.}, \bibinfo{author}{Eberhard, A.~C.},
  \bibinfo{author}{Engineer, F.~G.}, \bibinfo{author}{Fischetti, M.},
  \bibinfo{author}{Savelsbergh, M.~W.}, \& \bibinfo{author}{Tsoukalas, A.}
  (\bibinfo{year}{2014}).
\newblock \bibinfo{title}{Boosting the feasibility pump}.
\newblock {\it \bibinfo{journal}{Mathematical Programming Computation}\/},
  {\it \bibinfo{volume}{6}\/}, \bibinfo{pages}{255--279}.
\bibitem[{Boschetti et~al.(2009)Boschetti, Maniezzo, Roffilli \&
  R{\"o}hler}]{boschetti2009matheuristics}
\bibinfo{author}{Boschetti, M.~A.}, \bibinfo{author}{Maniezzo, V.},
  \bibinfo{author}{Roffilli, M.}, \& \bibinfo{author}{R{\"o}hler, A.~B.}
  (\bibinfo{year}{2009}).
\newblock \bibinfo{title}{Matheuristics: Optimization, simulation and control}.
\newblock In {\it \bibinfo{booktitle}{International Workshop on Hybrid
  Metaheuristics}\/} (pp. \bibinfo{pages}{171--177}).
\newblock \bibinfo{organization}{Springer}.
\bibitem[{Danna et~al.(2005)Danna, Rothberg \& Le~Pape}]{danna2005exploring}
\bibinfo{author}{Danna, E.}, \bibinfo{author}{Rothberg, E.}, \&
  \bibinfo{author}{Le~Pape, C.} (\bibinfo{year}{2005}).
\newblock \bibinfo{title}{Exploring relaxation induced neighborhoods to improve
  mip solutions}.
\newblock {\it \bibinfo{journal}{Mathematical Programming}\/},  {\it
  \bibinfo{volume}{102}\/}, \bibinfo{pages}{71--90}.
\bibitem[{Ehrgott(2005)}]{ehrgott2005multicriteria}
\bibinfo{author}{Ehrgott, M.} (\bibinfo{year}{2005}).
\newblock {\it \bibinfo{title}{Multicriteria optimization}\/} volume
  \bibinfo{volume}{491}.
\newblock \bibinfo{publisher}{Springer Science \& Business Media}.
\bibitem[{Ehrgott \& Gandibleux(2007)}]{ehrgott2007bound}
\bibinfo{author}{Ehrgott, M.}, \& \bibinfo{author}{Gandibleux, X.}
  (\bibinfo{year}{2007}).
\newblock \bibinfo{title}{Bound sets for biobjective combinatorial optimization
  problems}.
\newblock {\it \bibinfo{journal}{Computers \& Operations Research}\/},  {\it
  \bibinfo{volume}{34}\/}, \bibinfo{pages}{2674--2694}.
\bibitem[{Fernandes et~al.(2021)Fernandes, Goldbarg, Maia \&
  Goldbarg}]{fernandes2021multi}
\bibinfo{author}{Fernandes, I.~F.}, \bibinfo{author}{Goldbarg, E.~F.},
  \bibinfo{author}{Maia, S.~M.}, \& \bibinfo{author}{Goldbarg, M.~C.}
  (\bibinfo{year}{2021}).
\newblock \bibinfo{title}{Multi-and many-objective path-relinking: A taxonomy
  and decomposition approach}.
\newblock {\it \bibinfo{journal}{Computers \& Operations Research}\/},  (p.
  \bibinfo{pages}{105370}).
\bibitem[{Fischetti et~al.(2005)Fischetti, Glover \&
  Lodi}]{fischetti2005feasibility}
\bibinfo{author}{Fischetti, M.}, \bibinfo{author}{Glover, F.}, \&
  \bibinfo{author}{Lodi, A.} (\bibinfo{year}{2005}).
\newblock \bibinfo{title}{The feasibility pump}.
\newblock {\it \bibinfo{journal}{Mathematical Programming}\/},  {\it
  \bibinfo{volume}{104}\/}, \bibinfo{pages}{91--104}.
\bibitem[{Fischetti \& Lodi(2003)}]{fischetti2003local}
\bibinfo{author}{Fischetti, M.}, \& \bibinfo{author}{Lodi, A.}
  (\bibinfo{year}{2003}).
\newblock \bibinfo{title}{Local branching}.
\newblock {\it \bibinfo{journal}{Mathematical Programming}\/},  {\it
  \bibinfo{volume}{98}\/}, \bibinfo{pages}{23--47}.
\bibitem[{Fischetti \& Salvagnin(2009)}]{fischetti2009feasibility}
\bibinfo{author}{Fischetti, M.}, \& \bibinfo{author}{Salvagnin, D.}
  (\bibinfo{year}{2009}).
\newblock \bibinfo{title}{Feasibility pump 2.0}.
\newblock {\it \bibinfo{journal}{Mathematical Programming Computation}\/},
  {\it \bibinfo{volume}{1}\/}, \bibinfo{pages}{201--222}.
\bibitem[{Fonseca et~al.(2006)Fonseca, Paquete \&
  L{\'o}pez-Ib{\'a}nez}]{fonseca2006improved}
\bibinfo{author}{Fonseca, C.~M.}, \bibinfo{author}{Paquete, L.}, \&
  \bibinfo{author}{L{\'o}pez-Ib{\'a}nez, M.} (\bibinfo{year}{2006}).
\newblock \bibinfo{title}{An improved dimension-sweep algorithm for the
  hypervolume indicator}.
\newblock In {\it \bibinfo{booktitle}{2006 IEEE International Conference on
  Evolutionary Computation}\/} (pp. \bibinfo{pages}{1157--1163}).
\newblock \bibinfo{organization}{IEEE}.
\bibitem[{Forget et~al.(2022)Forget, Gadegaard \& Nielsen}]{forget2022warm}
\bibinfo{author}{Forget, N.}, \bibinfo{author}{Gadegaard, S.~L.}, \&
  \bibinfo{author}{Nielsen, L.~R.} (\bibinfo{year}{2022}).
\newblock \bibinfo{title}{Warm-starting lower bound set computations for
  branch-and-bound algorithms for multi objective integer linear programs}.
\newblock {\it \bibinfo{journal}{European Journal of Operational Research}\/},
  {\it \bibinfo{volume}{online first}\/}.
\bibitem[{Forget et~al.(2020)Forget, Klamroth, Gadegaard, Przybylski \&
  Nielsen}]{forget2020branch}
\bibinfo{author}{Forget, N.}, \bibinfo{author}{Klamroth, K.},
  \bibinfo{author}{Gadegaard, S.}, \bibinfo{author}{Przybylski, A.}, \&
  \bibinfo{author}{Nielsen, L.} (\bibinfo{year}{2020}).
\newblock \bibinfo{title}{Branch-and-bound and objective branching with three
  objectives}.
\newblock {\it \bibinfo{journal}{Preprint. Dec}\/}, .
\bibitem[{Gadegaard et~al.(2019)Gadegaard, Nielsen \&
  Ehrgott}]{gadegaard2019bi}
\bibinfo{author}{Gadegaard, S.~L.}, \bibinfo{author}{Nielsen, L.~R.}, \&
  \bibinfo{author}{Ehrgott, M.} (\bibinfo{year}{2019}).
\newblock \bibinfo{title}{Bi-objective branch-and-cut algorithms based on lp
  relaxation and bound sets}.
\newblock {\it \bibinfo{journal}{INFORMS Journal on Computing}\/},  {\it
  \bibinfo{volume}{31}\/}, \bibinfo{pages}{790--804}.
\bibitem[{Gandibleux et~al.(2021)Gandibleux, Gasnier \&
  Hanafi}]{gandibleuxprimal}
\bibinfo{author}{Gandibleux, X.}, \bibinfo{author}{Gasnier, G.}, \&
  \bibinfo{author}{Hanafi, S.} (\bibinfo{year}{2021}).
\newblock \bibinfo{title}{A primal heuristic to compute an upper bound set for
  multi-objective 0-1 linear optimisation problems}.
\newblock In \bibinfo{editor}{C.~F. Hayes}, \bibinfo{editor}{P.~Mannion}, \&
  \bibinfo{editor}{P.~Vamplew} (Eds.), {\it \bibinfo{booktitle}{Proc. of the
  1st Multi-Objective Decision Making Workshop (MODeM 2021)}\/}.
\bibitem[{Gandibleux et~al.(2003)Gandibleux, Morita \&
  Katoh}]{gandibleux2003impact}
\bibinfo{author}{Gandibleux, X.}, \bibinfo{author}{Morita, H.}, \&
  \bibinfo{author}{Katoh, N.} (\bibinfo{year}{2003}).
\newblock \bibinfo{title}{Impact of clusters, path-relinking and mutation
  operators on the heuristic using a genetic heritage for solving assignment
  problems with two objectives}.
\newblock In {\it \bibinfo{booktitle}{Proceedings of The Fifth Metaheuristics
  International Conference MIC’03}\/}.
\bibitem[{Gleixner et~al.(2021)Gleixner, Hendel, Gamrath, Achterberg, Bastubbe,
  Berthold, Christophel, Jarck, Koch, Linderoth et~al.}]{miplib}
\bibinfo{author}{Gleixner, A.}, \bibinfo{author}{Hendel, G.},
  \bibinfo{author}{Gamrath, G.}, \bibinfo{author}{Achterberg, T.},
  \bibinfo{author}{Bastubbe, M.}, \bibinfo{author}{Berthold, T.},
  \bibinfo{author}{Christophel, P.}, \bibinfo{author}{Jarck, K.},
  \bibinfo{author}{Koch, T.}, \bibinfo{author}{Linderoth, J.} et~al.
  (\bibinfo{year}{2021}).
\newblock \bibinfo{title}{{MIPLIB} 2017: data-driven compilation of the 6th
  mixed-integer programming library}.
\newblock {\it \bibinfo{journal}{Mathematical Programming Computation}\/},
  {\it \bibinfo{volume}{13}\/}, \bibinfo{pages}{443--490}.
\bibitem[{Glover(1997)}]{glover1997tabu}
\bibinfo{author}{Glover, F.} (\bibinfo{year}{1997}).
\newblock \bibinfo{title}{Tabu search and adaptive memory
  programming—advances, applications and challenges}.
\newblock In {\it \bibinfo{booktitle}{Interfaces in computer science and
  operations research}\/} (pp. \bibinfo{pages}{1--75}).
\newblock \bibinfo{publisher}{Springer}.
\bibitem[{Haimes(1971)}]{haimes1971bicriterion}
\bibinfo{author}{Haimes, Y.} (\bibinfo{year}{1971}).
\newblock \bibinfo{title}{On a bicriterion formulation of the problems of
  integrated system identification and system optimization}.
\newblock {\it \bibinfo{journal}{IEEE Transactions on Systems, Man, and
  Cybernetics}\/},  {\it \bibinfo{volume}{1}\/}, \bibinfo{pages}{296--297}.
\bibitem[{Hansen \& Mladenovi{\'c}(2001)}]{hansen2001variable}
\bibinfo{author}{Hansen, P.}, \& \bibinfo{author}{Mladenovi{\'c}, N.}
  (\bibinfo{year}{2001}).
\newblock \bibinfo{title}{Variable neighborhood search: Principles and
  applications}.
\newblock {\it \bibinfo{journal}{European Journal of Operational Research}\/},
  {\it \bibinfo{volume}{130}\/}, \bibinfo{pages}{449--467}.
\bibitem[{Kirlik \& Say{\i}n(2014)}]{kirlik2014new}
\bibinfo{author}{Kirlik, G.}, \& \bibinfo{author}{Say{\i}n, S.}
  (\bibinfo{year}{2014}).
\newblock \bibinfo{title}{A new algorithm for generating all nondominated
  solutions of multiobjective discrete optimization problems}.
\newblock {\it \bibinfo{journal}{European Journal of Operational Research}\/},
  {\it \bibinfo{volume}{232}\/}, \bibinfo{pages}{479--488}.
\bibitem[{Kiziltan \& Yucao{\u{g}}lu(1983)}]{kiziltan1983algorithm}
\bibinfo{author}{Kiziltan, G.}, \& \bibinfo{author}{Yucao{\u{g}}lu, E.}
  (\bibinfo{year}{1983}).
\newblock \bibinfo{title}{An algorithm for multiobjective zero-one linear
  programming}.
\newblock {\it \bibinfo{journal}{Management Science}\/},  {\it
  \bibinfo{volume}{29}\/}, \bibinfo{pages}{1444--1453}.
\bibitem[{Leitner et~al.(2016)Leitner, Ljubi{\'c}, Sinnl \&
  Werner}]{leitner2016ilp}
\bibinfo{author}{Leitner, M.}, \bibinfo{author}{Ljubi{\'c}, I.},
  \bibinfo{author}{Sinnl, M.}, \& \bibinfo{author}{Werner, A.}
  (\bibinfo{year}{2016}).
\newblock \bibinfo{title}{{ILP} heuristics and a new exact method for
  bi-objective 0/1 {ILPs}: Application to {FTTx}-network design}.
\newblock {\it \bibinfo{journal}{Computers \& Operations Research}\/},  {\it
  \bibinfo{volume}{72}\/}, \bibinfo{pages}{128--146}.
\bibitem[{L{\"o}hne \& Wei{\ss}ing(2017)}]{lohne2017vector}
\bibinfo{author}{L{\"o}hne, A.}, \& \bibinfo{author}{Wei{\ss}ing, B.}
  (\bibinfo{year}{2017}).
\newblock \bibinfo{title}{The vector linear program solver bensolve--notes on
  theoretical background}.
\newblock {\it \bibinfo{journal}{European Journal of Operational Research}\/},
  {\it \bibinfo{volume}{260}\/}, \bibinfo{pages}{807--813}.
\bibitem[{Pal \& Charkhgard(2019{\natexlab{a}})}]{pal2019feasibility}
\bibinfo{author}{Pal, A.}, \& \bibinfo{author}{Charkhgard, H.}
  (\bibinfo{year}{2019}{\natexlab{a}}).
\newblock \bibinfo{title}{A feasibility pump and local search based heuristic
  for bi-objective pure integer linear programming}.
\newblock {\it \bibinfo{journal}{INFORMS Journal on Computing}\/},  {\it
  \bibinfo{volume}{31}\/}, \bibinfo{pages}{115--133}.
\bibitem[{Pal \& Charkhgard(2019{\natexlab{b}})}]{pal2019fpbh}
\bibinfo{author}{Pal, A.}, \& \bibinfo{author}{Charkhgard, H.}
  (\bibinfo{year}{2019}{\natexlab{b}}).
\newblock \bibinfo{title}{{FPBH}: A feasibility pump based heuristic for
  multi-objective mixed integer linear programming}.
\newblock {\it \bibinfo{journal}{Computers \& Operations Research}\/},  {\it
  \bibinfo{volume}{112}\/}, \bibinfo{pages}{104760}.
\bibitem[{Parragh et~al.(2009)Parragh, Doerner, Hartl \&
  Gandibleux}]{parragh2009heuristic}
\bibinfo{author}{Parragh, S.~N.}, \bibinfo{author}{Doerner, K.~F.},
  \bibinfo{author}{Hartl, R.~F.}, \& \bibinfo{author}{Gandibleux, X.}
  (\bibinfo{year}{2009}).
\newblock \bibinfo{title}{A heuristic two-phase solution approach for the
  multi-objective dial-a-ride problem}.
\newblock {\it \bibinfo{journal}{Networks: An International Journal}\/},  {\it
  \bibinfo{volume}{54}\/}, \bibinfo{pages}{227--242}.
\bibitem[{Requejo \& Santos(2017)}]{requejo2017feasibility}
\bibinfo{author}{Requejo, C.}, \& \bibinfo{author}{Santos, E.}
  (\bibinfo{year}{2017}).
\newblock \bibinfo{title}{A feasibility pump and a local branching heuristics
  for the weight-constrained minimum spanning tree problem}.
\newblock In {\it \bibinfo{booktitle}{International Conference on Computational
  Science and Its Applications}\/} (pp. \bibinfo{pages}{669--683}).
\newblock \bibinfo{organization}{Springer}.
\bibitem[{Soylu(2015)}]{soylu2015heuristic}
\bibinfo{author}{Soylu, B.} (\bibinfo{year}{2015}).
\newblock \bibinfo{title}{Heuristic approaches for biobjective mixed 0--1
  integer linear programming problems}.
\newblock {\it \bibinfo{journal}{European Journal of Operational Research}\/},
  {\it \bibinfo{volume}{245}\/}, \bibinfo{pages}{690--703}.
\bibitem[{Tamby \& Vanderpooten(2021)}]{tamby2021enumeration}
\bibinfo{author}{Tamby, S.}, \& \bibinfo{author}{Vanderpooten, D.}
  (\bibinfo{year}{2021}).
\newblock \bibinfo{title}{Enumeration of the nondominated set of multiobjective
  discrete optimization problems}.
\newblock {\it \bibinfo{journal}{INFORMS Journal on Computing}\/},  {\it
  \bibinfo{volume}{33}\/}, \bibinfo{pages}{72--85}.
\bibitem[{Zitzler \& Thiele(1999)}]{zitzler1999multiobjective}
\bibinfo{author}{Zitzler, E.}, \& \bibinfo{author}{Thiele, L.}
  (\bibinfo{year}{1999}).
\newblock \bibinfo{title}{Multiobjective evolutionary algorithms: a comparative
  case study and the strength pareto approach}.
\newblock {\it \bibinfo{journal}{IEEE Transactions on Evolutionary
  Computation}\/},  {\it \bibinfo{volume}{3}\/}, \bibinfo{pages}{257--271}.
\bibitem[{Zitzler et~al.(2003)Zitzler, Thiele, Laumanns, Fonseca \&
  Da~Fonseca}]{zitzler2003performance}
\bibinfo{author}{Zitzler, E.}, \bibinfo{author}{Thiele, L.},
  \bibinfo{author}{Laumanns, M.}, \bibinfo{author}{Fonseca, C.~M.}, \&
  \bibinfo{author}{Da~Fonseca, V.~G.} (\bibinfo{year}{2003}).
\newblock \bibinfo{title}{Performance assessment of multiobjective optimizers:
  An analysis and review}.
\newblock {\it \bibinfo{journal}{IEEE Transactions on Evolutionary
  Computation}\/},  {\it \bibinfo{volume}{7}\/}, \bibinfo{pages}{117--132}.

\end{thebibliography}

\end{document}